\documentclass[12pt]{article}
\usepackage{amsfonts}
\usepackage{amsmath}
\usepackage{epsfig}

\title{Symplectic Tiling Billiards, Planar Linkages, and Hyperbolic
  Geometry}
\author{Richard Evan Schwartz \thanks{Supported by N.S.F. Grant DMS-2102802}}

\newtheorem{theorem}{Theorem}[section]

\newtheorem{lemma}[theorem]{Lemma}

\newtheorem{conjecture}[theorem]{Conjecture}

\def\startproof{{\bf {\medskip}{\noindent}Proof: }}

\def\endproof{$\spadesuit$  \newline}

\def\C{\mbox{\boldmath{$C$}}}%
\def\H{\mbox{\boldmath{$H$}}}%
\def\Q{\mbox{\boldmath{$Q$}}}%
\def\R{\mbox{\boldmath{$R$}}}%
\def\Z{\mbox{\boldmath{$Z$}}}%

\begin{document}
\maketitle

\begin{abstract}
 In this paper I will unite two games,
  symplectic billiards and tiling billiards. The new
  game is called symplectic tiling billiards.
  I will prove a result about periodic orbits of
  symplectic tiling billiards in a
  very special case and then show how this result
  combines with the construction in 
  Thurston's paper {\it Shapes of Polyhedra\/} to
  give hyperbolic structures on moduli spaces
  of planar equilateral polygons.  One corollary
  is that the configuration space of the hexagonal
  planar linkage with unit-length rods (modulo
  isometry)  has an
  algebraically defined hyperbolic structure in which
  it is a $10$-cusped hyperbolic $3$-manifold
  that is tiled by $15$ regular ideal octahedra.
  The $10$ cusps correspond to the $10$ maximally degenerate
  configurations.
      \end{abstract}

\section{Introduction}

Billiards, of course, needs no introduction.
However, it has two exotic cousins which are less
well known, {\it symplectic billiards\/} and
{\it tiling billiards\/}.   In this paper I will
unite these two topics. I
call the new game {\it symplectic tiling billiards\/}.
Perhaps anyone who knows about both
symplectic billiards and tiling billiards could stop reading
now and define symplectic tiling billiards for themselves just
based on the name.

For ease of exposition
I will stick to the polygonal cases of all these
topics.
Symplectic billiards is perhaps best played on
a pair of polygons, $A$ and $B$, as shown in
Figure 1.1.  Starting with a pair
$(a_1,b_2) \in \partial A \times \partial B$ one
produces a pair $(a_3,b_4) \in \partial A \times \partial B$
using the rule below.

\begin{center}
\resizebox{!}{1.3in}{\includegraphics{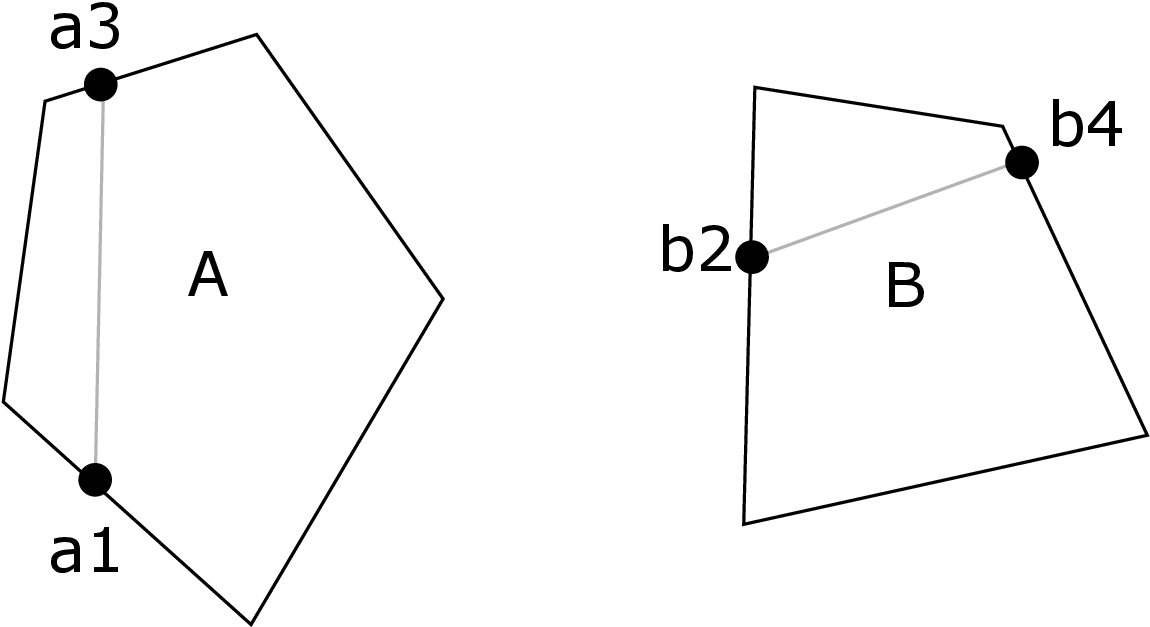}}
\newline
{\bf Figure 1.1:\/} Symplectic Billiards Defined
\end{center}

In words, the line connecting
$a_1$ to $a_3$ is parallel to the side of $B$
containing $b_2$ and the line connecting
$b_2$ to $b_4$ is parallel to the side of $A$
containing $a_3$.   One then iterates and
considers the dynamics.   I first learned about
symplectic billiards from Peter Albers and
Serge Tabachnikov.  We 
later wrote a paper [{\bf ABSST\/}] about the
subject, proving a few foundational results.
The two-table perspective is explored extensively in the more recent
work [{\bf ALW\/}].

Tiling billiards is a variant of billiards played
on the edges of a planar tiling.  Figure 1.2 shows
the rule. 

\begin{center}
\resizebox{!}{2.2in}{\includegraphics{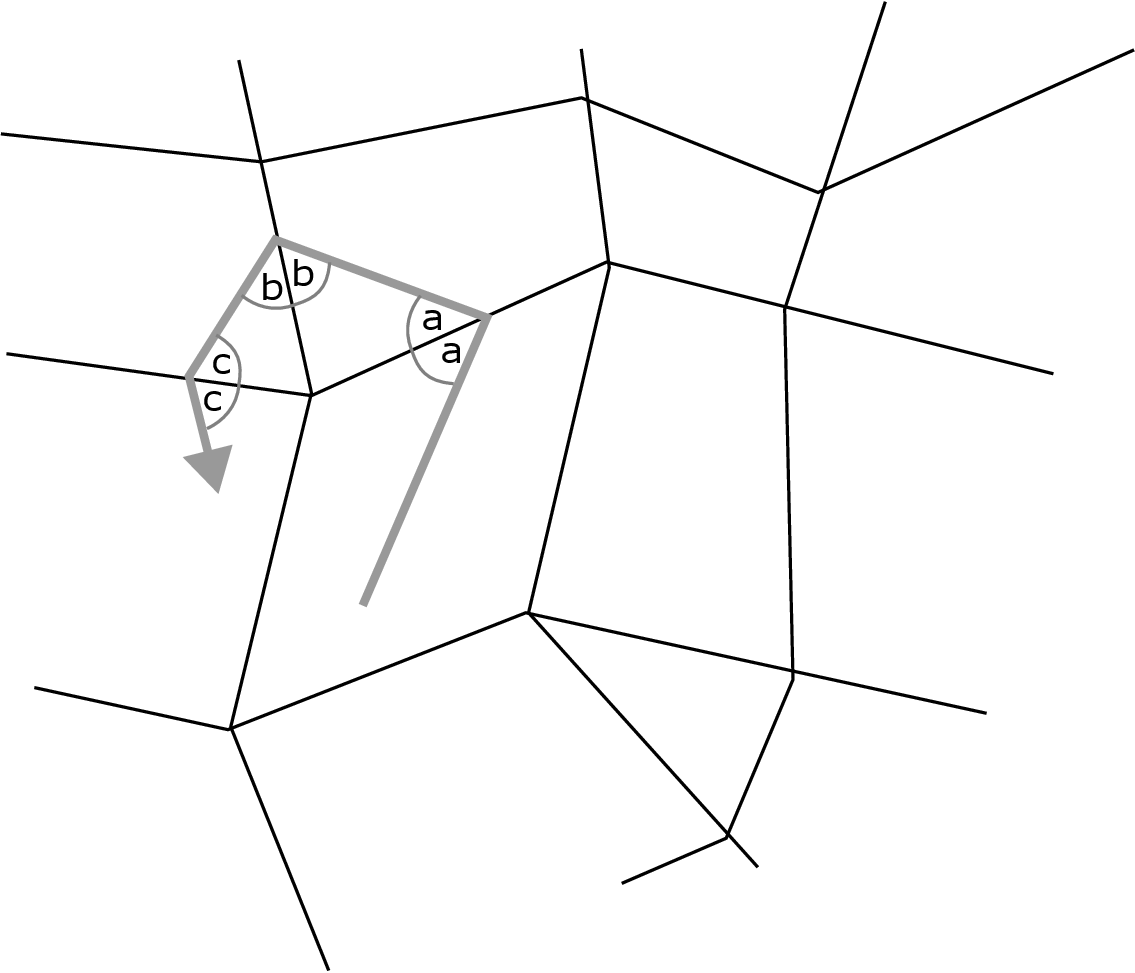}}
\newline
{\bf Figure 1.2:\/} Tiling Billiards
\end{center}

The rule is essentially the same as for
billiards, except that the trajectory
refracts through the edges rather than
bouncing off them.  I first learned about
tiling billiards from Serge Tabachnikov.
Now there is a growing literature on the
subject.  See
[{\bf BDFI\/}] and the references therein.

In \S 2 I will define {\it symplectic tiling
billiards\/} and make a few general
remarks about it. I will also show
the results of a few easy experiments.
The game is played relative to a pair of
tilings of the plane, though one could
specialize to the case where the two
tilings are the same.

In \S 3 I will consider a special case of
this game, a kind of ``local version'',
in which the planar tilings
involved each consist of $N$
infinite sectors bounded by $N$ rays emanating
from the origin.  We call such a tiling an
$N$-{\it sunburst\/}.   We call the $N$-sunburst
{\it balanced\/} if the sum of the $N$  unit vectors
parallel to the rays is $0$.   As a special
case, we call the $N$-sunburst {\it regular\/}
if the rays are parallel to the $N$th roots of unity.
Figure 1.3 shows a pair $(A,B)$  of $7$-sunbursts
where $A$ is regular and $B$ is balanced. The rays
of $A$ point outwards and the rays of $B$ point inwards,
as indicated by the arrows.

\begin{center}
\resizebox{!}{2.2in}{\includegraphics{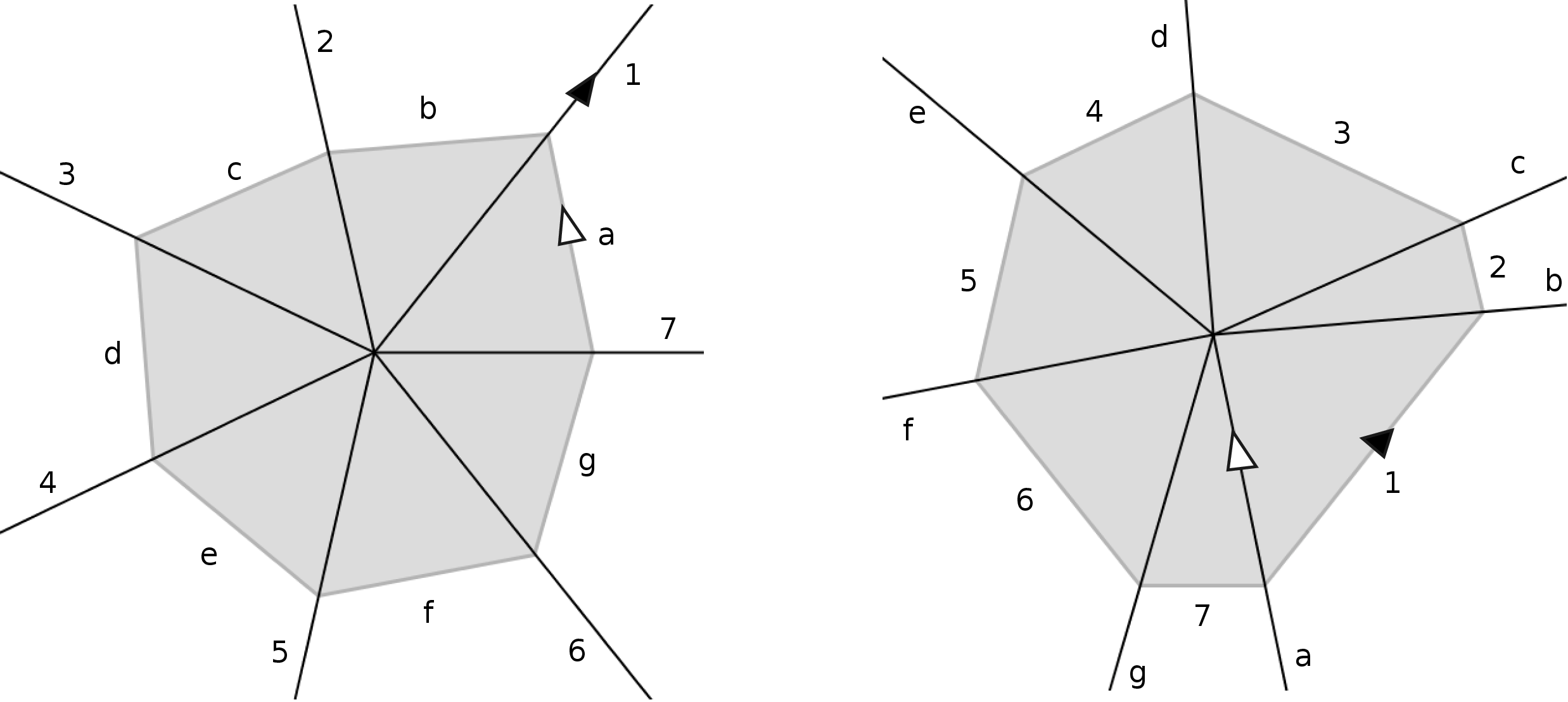}}
\newline
{\bf Figure 1.3:\/} Adapted polygons supported by a pair of
$7$-sunbursts
\end{center}

Figure 1.3 shows a very special situation.  Each sunburst
has an inscribed counterclockwise-oriented convex polygon such that the edges and
rays having the same label are parallel in the oriented sense
indicated by the arrows.  Let us provisionally
call such polygons {\it adapted polygons\/} and 
say that the pair of sunbursts {\it supports\/} them.
Later on
we will see that these polygons are  periodic orbits for
symplectic tiling billiards.

We say that a pair $(A,B')$ of sunbursts is a
{\it phase modification\/} of the pair $(A,B)$ if
$B'$ is obtained from $B$ by rotating about the origin.
Here is a restatement of one of our main results,
Theorem \ref{two}.

\begin{theorem}
  \label{local}
  A pair $(A,B)$ of $N$-sunbursts, with $A$ regular and $B$ balanced,
  has a unique phase modification which supports adapted $N$-gons.
\end{theorem}

In \S 4 we consider planar equilateral polygons.
All the edges have the same length.
If we fix the number $N$ of sides, then the
moduli space of equilateral $N$-gons modulo
similarity is equivalent to the configuration
space of the mechanical linkage made from
$N$ unit-length rods.
Theorem \ref{local} gives a very clean
bijection between similarity
classes of strictly convex equilateral $N$-gons and
similarity classes of strictly convex equiangular $N$-gons.

Here is the idea.   Let $L$ be an equilateral $N$-gon.
By taking the rays parallel to the
edges of $L$, in order, we get a balanced
$N$-sunburst $B$.  We let $A$ be the regular
$N$-sunburst.  We apply Theorem \ref{local}
to $(A,B)$ to get a phase modification $(A,B')$ which
supports
adapted polygons.  The similarity class
of the adapted $N$-gon inscribed in $B'$ (as in the right side of
Figure 1.3) gives
us our class of equiangular $N$-gon $P_L$.
Our association respects the equivalence class
and this gives us our bijection $[L] \to [P_L]$.

Our correspondence
is akin to the one given by
Misha Kapovich and John Millson [{\bf KM\/}], but it
is more direct, more algebraic, and easier to compute.
See \S \ref{compute} and \S \ref{algXX}.
The Kapovich-Millson
corrspondence involves the Riemann mapping
theorem, a transcendental construction.

In \S \ref{hypXX}  we recall some features of
William Thurston's famous {\it Shapes of Polyhedra\/}
construction [{\bf T\/}] of a complex hyperbolic structure
on the moduli space of polyhedra with prescribed
cone angles.
When we restrict our attention to
polyhedra which are doubles of strictly
convex equiangular polygons, we get a
real hyperbolic structure on the space of
similarity classes of
strictly convex equiangular polygons.
All this is well-known,
and I will give a self-contained account
in \S \ref{hypXX}.
Thurston's construction realizes the moduli space
of strictly convex equiangular $N$-gons as the
interior of convex polytope in hyperbolic space $\H^{N-3}$.
The polytope is canonically defined up to the action
of algebraic matrices. Hence it makes
sense to talk about algebraic points in this polytope.
 When $N$ is odd the polytope
is bounded and when $N$ is even it has some
ideal vertices.

We use our correspondence
$[L]\to [P_L]$ to give
a hyperbolic structure to the moduli
space ${\cal C\/}_N$
of strictly convex equilateral $N$-gons.
We call this the {\it algebraic hyperbolic structure\/}.
Equipped with the algbraic hyperbolic structure,
${\cal C\/}_N$ is just the Thurston polytope in
$\H^{N-3}$.
Our correspondence just imports the
Thurston construction to the equilateral case.

\begin{theorem}
  \label{algebra}
  Relative to the algebraic hyperbolic structure,
  a similarity class in ${\cal C\/}_N$ has a
  representative  with algebraic coordinates
    if and only if the class has algebraic coordinates
  in $\H^{N-3}$.
\end{theorem}

Our construction only
works in the strictly convex case but there is
trick to extend the construction to the general case.
Using the action of the permutation group,
which acts on the space ${\cal A\/}_N$
of all equilateral $M$-gons,
we can extend our hyperbolic structure on
${\cal C\/}_N$ to one on ${\cal A\/}_N$.
See \S \ref{put}.
In the even case, we need to adjoin the
ideal vertices to ${\cal A\/}_N$ in order to
include the degenerate polygons that lie in a
single line.

When $N>6$ the space ${\cal A\/}_N$ has
various conical singularities because the
various copies of ${\cal C\/}_N$ do not
fit nicely together around codimension $2$
faces.
The picture of ${\cal A\/}_N$ is very satisfying
when $N=5,6$.

\begin{theorem}
  \label{penta}
  Relative to the algebraic hyperbolic structure,
    ${\cal A\/}_5$  is a
hyperbolic surface of Euler characteristic $-3$.
The space is tiled by $12$ regular right angled
pentagons which meet $4$ around
each vertex.
\end{theorem}

\begin{theorem}
  \label{hex}
  Relative to the algebraic hyperbolic structure,
  ${\cal A\/}_6$  is a finite volume
  $10$-cusped hyperbolic $3$-manifold
  that is tiled by $15$ regular ideal
  octahedra which meet $4$ around each edge.
\end{theorem}

One thing that inspired me to
consider equilateral polygons is that I had recently
heard a great talk given by 
Juergen Richter-Gebert [{\bf R-G\/}]
about his hyperbolic structure on the space of
equilateral pentagons.
Richter-Gebert has a different way to
give a hyperbolic structure in the pentagonal case.
His construction, like that in [{\bf KM\/}], is
transcendental and seemingly hard to compute.

I would describe the original version of the
paper as a meal that I threw
together based on ideas that were dropped
on my plate while I dined in Heidelberg
and Marseille during a very happy summer
in 2023.  I probably had the key idea for this paper
while in free-fall
riding the Hurricane Loop waterslide at
Miramar water park in Weinheim.
Thanks to the probing comments of
the anonymous referees, and also thanks to a
great insight of Jannik Westermann
which I describe in \S \ref{leftright},
this version of the paper is both deeper
and sharper than the original.

I thank Peter Albers, Diana Davis, Peter Doyle,
Aaron Fenyes, Fabian Lander, 
Juergen Richter-Gebert, Joe Silverman,
Sergei Tabachnikov,
Jannik Westermann, and two anonymous referees
for helpful discussions about this paper.
Finally, I am grateful for the support I've had
during this time period from the University of Heidelberg,
from CIRM (Luminy), from
the National Science Foundation, and from the Simons
Foundation.

\newpage

\section{Symplectic Tiling Billiards}

\subsection{Basic Definition}

For us, a {\it tiling\/} is a subdivision of the plane into
convex polygonal regions.  These polygonal regions
are allowed to be unbounded.    The simplest unbounded
case is that of a {\it sector\/}, namely a region bounded
by rays which make an angle
of less than $\pi$ with each other.

Given a tiling $A$, a {\it particle\/} on $A$ is a
point contained on the interior of an edge of $A$,
together with a choice of a direction pointing into
one of the two regions of $A$ adjacent to $e$.
One could encode the direction by a vector
transverse to $e$.

We say that two tilings $A$ and $B$ are
{\it transverse\/} if no edge of $A$ is parallel to
an edge of $B$.
Suppose that $(A,B)$ are transverse and
$(a_1,b_2)$ are a pair of particles, with $a_1$ being
a particle of $A$ and $b_2$ being a particle of $B$.

We define $a_3$ as follows.  The line connecting
$a_1,a_3$ is parallel to the edge of $B$ containing
$b_2$.  The direction at $a_3$ goes in the same
direction as the direction at $a_1$.  That is, one
and the same vector along the line $\overline{a_1a_3}$
would serve as a transverse vector.   We define
$b_4$ in the same way, swapping the roles of
$A$ and $B$.  Figure 2.1 shows the construction.

\begin{center}
\resizebox{!}{2.8in}{\includegraphics{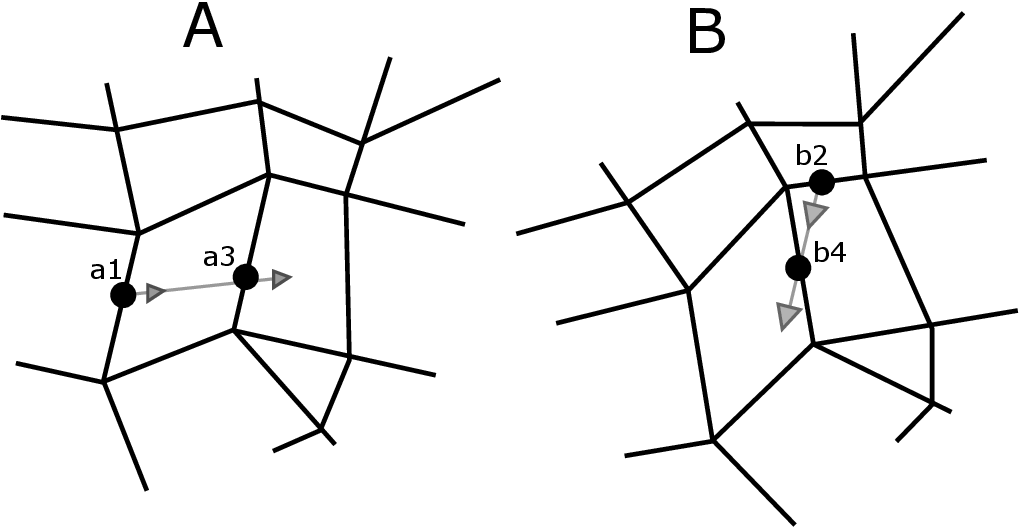}}
\newline
{\bf Figure 2.1:\/} Symplectic Billiards Defined
\end{center}

\subsection{Remarks on the Definition}

Here are some comments about the basic definition.
\newline
\newline
\noindent
{\bf Unbounded Tiles:\/}
One subtle point about this definition is that
the point $a_3$ or the point $b_4$ might not be
defined in case $A$ or $B$ has unbounded tiles.
What happens here is that the relevant ray
simply heads off to infinity without intersecting
an edge of the tiling.  We allow this, and indeed
it might present an interesting case to study, but the
squeamish reader could avoid this problem by
only considering tilings with bounded tiles.
\newline
\newline
{\bf Lack of Transversality:\/}
One might also want to consider the case when
$A$ and $B$ are not transverse.  For instance,
one might like to play this game on a single tiling,
setting $A=B$.  As in symplectic billiards, one
requires that the particles
$a_1$ and $b_2$ are not contained in parallel edges.
\newline
\newline
{\bf Affine Symmetry:\/}
Like symplectic billiards, symplectic tiling
billiards is affinely natural.  If $T$ is an
affine transformation of the plane, then
$T$ maps the orbits relative to the pair
$(A,B)$ to the orbits relative to the pair
$(A',B')$ where $A'=T(A)$ and $B'=T(B)$.
Also, if $T$ is a dilation
then the orbits relative to $(A,B)$ are the
same as the orbits relative to $(A,B')$.

It might be interesting to study symplectic
tiling billiards on tilings which have affine
symmetry.  These are called
affine crystallographic groups.
\newline
\newline
{\bf Half Translation Surfaces:\/}
Symplectic tiling billiards can also be
played on a torus.   This is equivalent to
considering the game relative to a pair
of doubly periodic tilings, and then
considering the orbits on the quotient space.

More generally, one can play the game relative
to a half-translation surface.  Recall that a half-translation
surface is a metric on a surface in which all but
a discrete set of points are locally isometric
to the Euclidean plane and the remaining
points are cone points having cone angle
$\pi k$ for various integers $k$.
One additional requirement for these surfaces
is that there is a global parallel line field.
(This is not quite implied by the other
conditions.)

Let's say that a tiling of a translation surface
is a decomposition of the surface into convex
polygons such that every cone point appears
as a vertex.   Other points might be vertices
as well.  Choices of globally parallel line
fields would give a way to line up the two
tilings.

\subsection{Rotated Square Grids}

Here I make some remarks about some
experiments I did with
symplectic tiling billiards in
the case when $A$ and $B$ are both
square grids.  In this case, the only
parameter is the way $A$ and $B$ are
rotated with respect to each other.

Given
$t \in \Q$ define
\begin{equation}
  \label{RAT}
  z_t=\frac{1-t^2}{1+t^2} + i \bigg(\frac{2t}{1+t^2}\bigg).
\end{equation}
This is the usual rational parametrization of the unit circle.
We normalize so that $A$ is the usual
square grid.
Let $A_t$ denote the result of multiplying the usual
square grid by $z_t$.
The program I wrote uses exact rational
arithmetic to explore this case of
rationally rotated square grids.

For the pair $(A,A_0)$, which is the same
as just playing the game on $A$, the
orbits are just rays.
The choice $t=1/3$ yields $z_t=(4/5)+(3/5)i$.
This is the simplest non-trivial rational case.
It is based on the $(3,4,5)$ right triangle.
Let's take a look.

Figure 4 shows a picture of an orbit
with respect to $(A,A_{1/3})$.
The picture on
the right shows a close-up of the most
complicated part of the orbit on the right.
(I have also continued the orbit a bit further
on the right.)
The slightly thicker segments on the right
are actually unions of extremely close
parallel segments.

\begin{center}
\resizebox{!}{2.5in}{\includegraphics{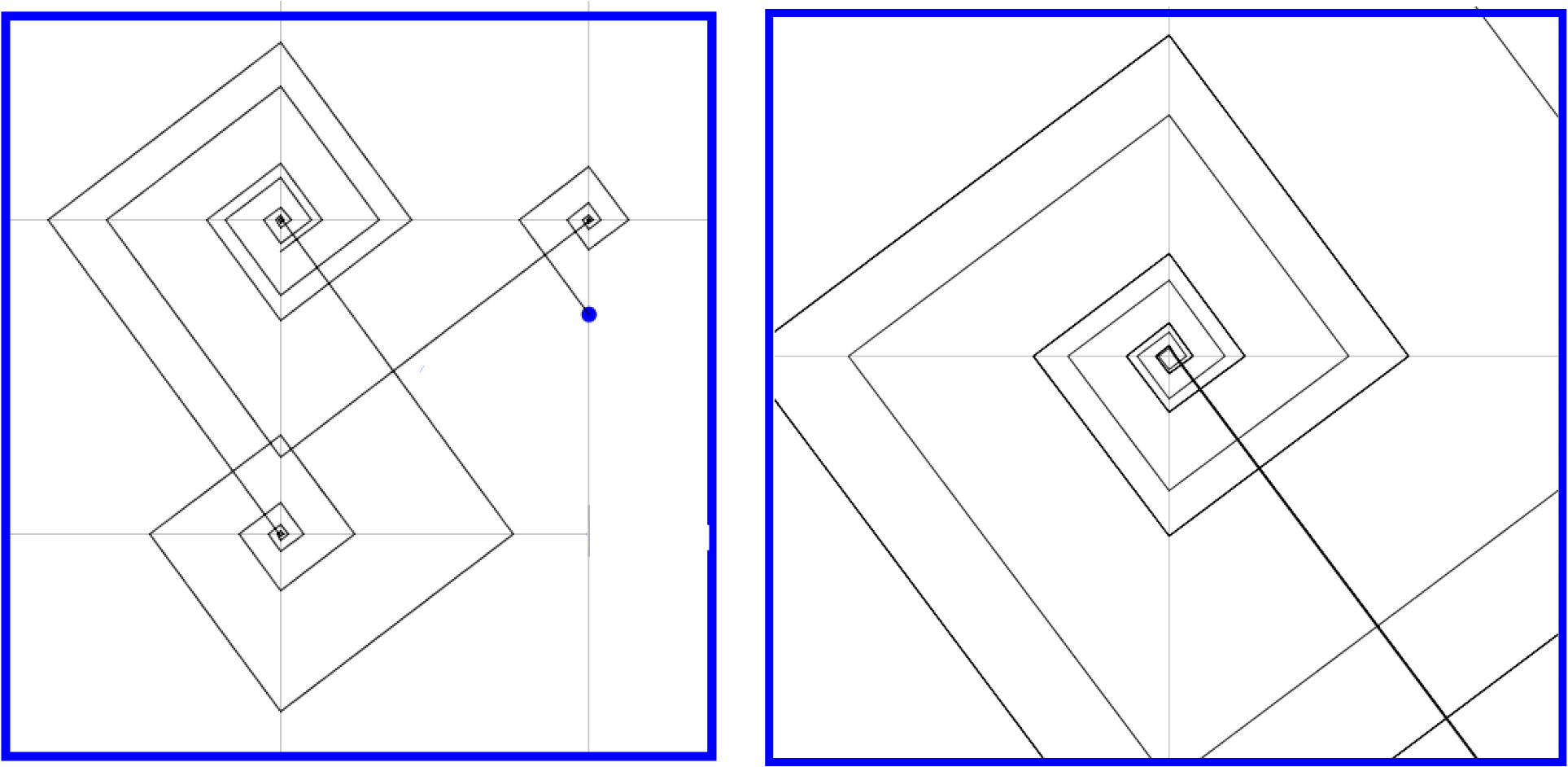}}
\newline
    {\bf Figure 2.2:\/} The left half of an orbit on $(A,A_{1/3})$.
\end{center}

The orbit seems to be bounded and aperiodic.  
I did not attempt to prove this, but
my rational calculations reveal that
the numerators and denominators of the
coordinates of the vertices are
tending to $\infty$.  Geometrically, the orbit has
an attracting limit cycle.  All this
would not be hard to prove; one would
look at the first return map to a suitable
interval, get an interval exchange
transformation, and check that it had
an attracting fixed point.  Here is a concrete
conjecture.

\begin{conjecture}
All orbits on $(A,A_{1/3})$ are bounded and
get attracted to a limit cycle.
The limit cycle itself is a periodic orbit.
\end{conjecture}

Theorem \ref{LR} gives some justification for the
nature of Figure 2.2 above and Figure 2.3 below.

\begin{lemma}
  \label{shoot}
  For a pair of rationally rotated grids, it is impossible
  for an orbit on the left (or on the right) to remain
  forever on the $4$ edges incident to a single vertex.
\end{lemma}

\startproof 
We will suppose that this happens on the left.
Let $V_L$ be the vertex on the left that
the orbit moves around.  As the orbit
on the left side moves around $V_L$ is
must always make a $\pi/2$ degree turn,
either clockwise or counterclockwise.
But then the orbit on the right must forever
hit horizontal and vertical edges in alternation.
But then there is a vertex $V_R$ on the right
such that the orbit on the right stays on the
$4$ edges incident to $V_R$.

Suppose the orbit on the left gets closer to
$V_L$ after $4$ steps.  Then, by symmetry,
the orbit spirals in towards $V_L$ in a
geometric series: the distance to $V_L$
drops by a definite factor $\lambda<1$ at
each revolution.   Theorem \ref{LR} 
says that on the right the distance from
the orbit to $V_R$ increases by $1/\lambda$
after each revolution.  (This result would
be easy to work out by hand in our setting here.)
But then the orbit
on the right eventually escapes and we
have a contradiction.

The other possibility is that the orbit on the
left is periodic.   But then the grids are
rotated by $\pi/4$ degrees relative to
each other.  These are not rationally
rotated grids.
\endproof

Lemma \ref{shoot} gives some explanation for
why the orbits in Figure 2.2 and Figure 2.3 seem
to spiral close to a vertex and then suddenly
shoot out.
\newline

Figure 2.3 shows a picture of an unbounded orbit
on $(A,A_{7/11})$.  I picked this parameter somewhat
randomly.  After making two big spirals, the orbit
starts heading southeast in a periodic pattern with
a drift.  The orbit is clearly unbounded, though I
did not attempt a proof.

\begin{center}
\resizebox{!}{4.5in}{\includegraphics{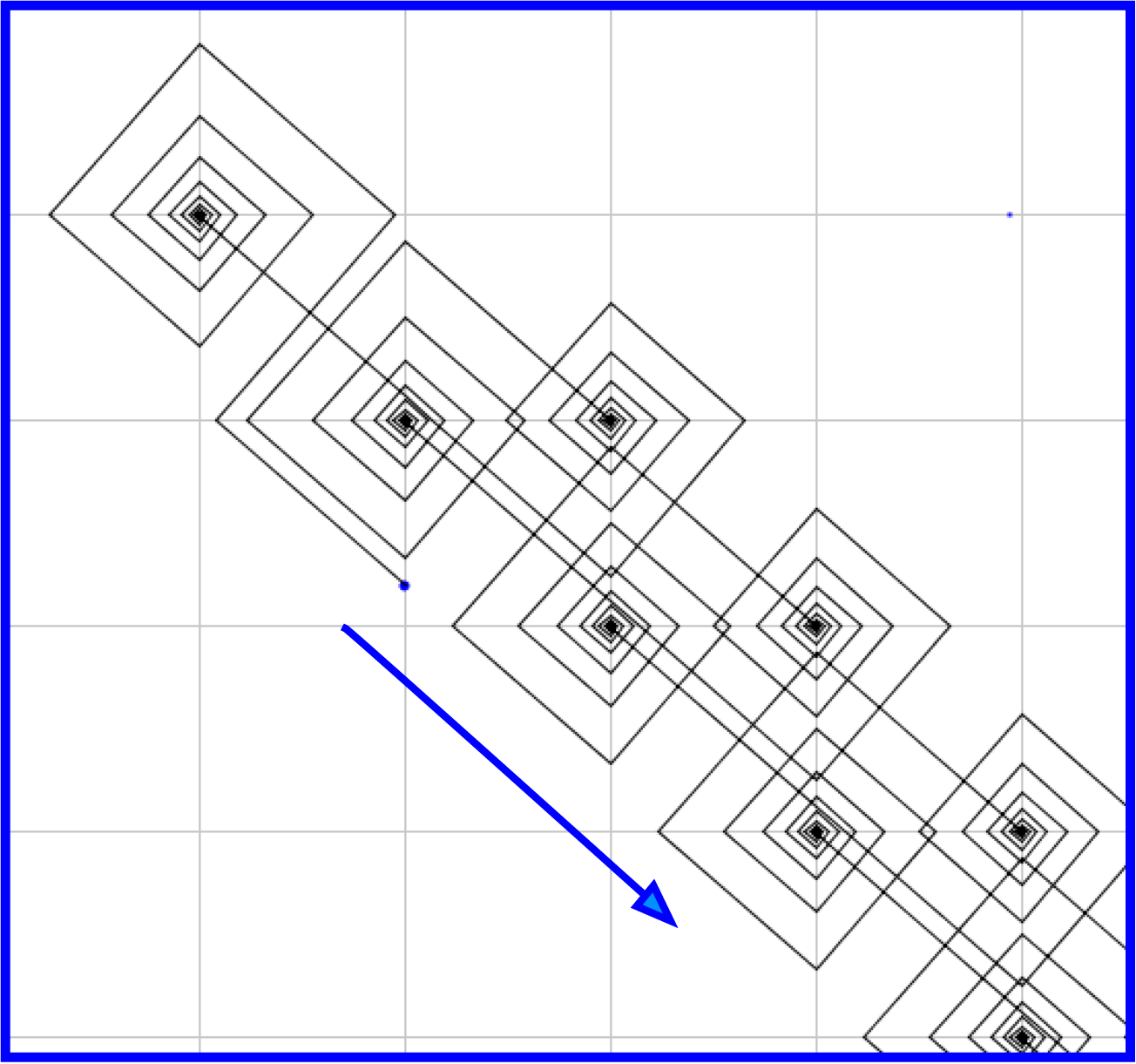}}
\newline
    {\bf Figure 2.3:\/} A symplectic billiard orbit on $(A,A_{7/11})$.
\end{center}

Some parameters seem to support both unbounded orbits
and bounded orbits. It would be nice to classify the
rational parameters according to which kinds of
orbits they support.

For irrationally rotated square grids,
one can sometimes get entirely periodic orbits.
For instance, if we play on $(A,B)$, where $B$ is obtained
by rotating $A$ $\pi/4$ radians, then all orbits are
periodic, and they make squares in each factor.
Compare the end of
the proof of Lemma \ref{shoot}.

\newpage

\section{Sunbursts}

\subsection{Basic Definitions}
\label{mainres}

The goal in this chapter is to prove
Theorem \ref{local} and related results.
An $N$-{\it sunburst\/} is a union of $N$ rays emanating from
the origin such that the convex hull of the rays is the whole plane.
An $N$-sunburst defines a tiling in the plane in which the
tiles are unbounded sectors based at the origin.
In this section we will consider symplectic tiling billiards
with respect to two $N$-sunbursts $A$ and $B$.
The number $N$ is the same for both $A$ and $B$.

We orient the rays of $A$ outward
and the rays of $B$ inward, as shown in Figure 3.1.
(Compare Figure 1.3.)
For the entire chapter, we restrict our attention to
the situation where
we have an orbit that starts with $a_1 \in A_1$ and $b_2 \in B_2$,
so that the particle at $a_1$ points into the sector bounded
by $A_1, A_3$ and the particle at $b_2$ points into the
sector bounded by $B_2, B_4$.  When we say that
$(A,B)$ has periodic orbits, we implicitly mean this kind.
By dilation symmetry, one orbit on $(A,B)$ is periodic
if and only if they all are.

\begin{center}
\resizebox{!}{2.6in}{\includegraphics{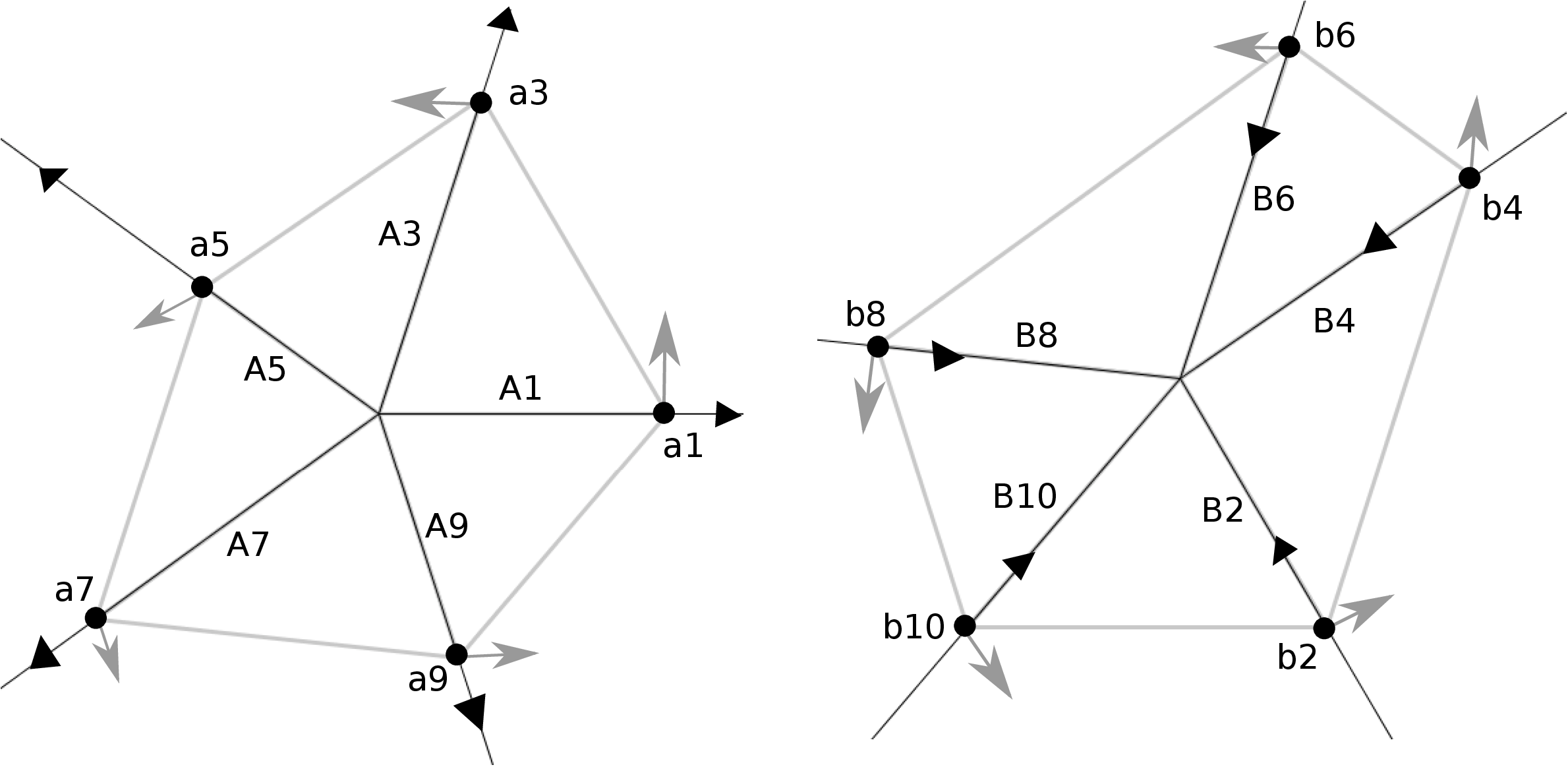}}
\newline
{\bf Figure 3.1:\/}  A periodic orbit relative to a pair of
$5$-sunbursts
\end{center}

Let $(A,B)$ be a pair of sunbursts.  Call an
orbit $\cal O$ of $(A,B)$ {\it woven\/} if the restriction
of ${\cal O\/}$ to $A$, which we call
${\cal O\/}_A$, circulates counterclockwise
around $A$ and
if (with the obvious notation) ${\cal O\/}_B$
circulates counterclockwise around $B$.
The orbit shown in Figure 3.1 is both
woven and periodic.  In this situation
${\cal O\/}_A$ and ${\cal O\/}_B$ are
both convex polygons.

\subsection{A Criterion for Periodicity}
\label{leftright}

Say that ${\cal O\/}$ is {\it left-convex\/}
(respectively {\it right-convex\/})
if ${\cal O\/}_A$ (respectively ${\cal O\/}_B$) is a closed convex polygon. 
In the first version of this paper, I considered
left-convex orbits but did not inquire as to
whether left-convex orbits were also
right-convex and hence periodic.
However, Jannik Westermann read the first
version of the paper and asked this question.
He  noticed that a woven orbit seems to be
left-convex if and only if it is right-convex.
Jannik gave a geometric proof of this
fact for pairs of $3$-sunbursts.  Subsequently, I found
a proof of the general case. 

    \begin{theorem}
      \label{LR}
      An woven orbit relative to a pair of $N$-sunbursts is left-convex if
      and only if it is right-convex if and only if it is periodic.
    \end{theorem}

    The rest of this section is devoted to proving Theorem \ref{LR}.
    Our proof goes
    through some elementary complex analysis.
    I also used this idea in [{\bf Sch2\/}].
    I'd prefer a geometric proof, but I don't have one.
        \newline
    
A {\it complex $N$-sunburst\/} is an
ordered list of $N$ complex lines through the origin
in $\C^2$.  We will usually drop the word {\it complex\/}
in our discussion. We will be interested in a pair
$(A,B)$ of $N$-sunbursts.  We write these
as $A=A_1,A_3,...,A_{2N-1}$ and
$B=B_2,B_4,...,B_{2N}$.

We say that $(A,B)$ is a {\it good pair\/} if
$A_i \not = B_{i \pm 1}$ for all indices.
Let $(A,B)$ be a good pair.
Given $z_1 \in A_1-\{0\}$ we let
$z_3= B_2' \cap A_3$ where $B_2'$ is
the complex line through $z_1$ parallel to $B_2$.
Since $(A,B)$ is a good pair, we have
$z_3 \in A_3-\{0\}$.  In the same way
we define $z_5=B_4' \cap A_5 $, etc.
This gives us points
$z_7,...,z_{2N-1},z_{2N+1}$.
The ratio
$$\lambda_A=z_{2N+1}/z_1$$
  makes sense because both points
  lie in $A_1$.  Also, by scaling symmetry,
  $\lambda_A$
  is independent of the choice of $z_1$.
  We would get the same value
  if we started with $z_3 \in A_3$ and set
  $\lambda_A=z_{2N+3}/z_3$. Etc.
    Likewise we define $\lambda_B$.
  
  \begin{theorem}
    \label{main}
    We have $\lambda_A \lambda_B=1$ for all good pairs $(A,B)$.
  \end{theorem}

  Theorem \ref{main} applies
  in the real case to the pairs of
  sunbursts considered in Theorem \ref{LR}.
  The corresponding woven orbits are left-convex if and only
  if $\lambda_A=1$ and right-convex if and only if
  $\lambda_B=1$.  But Theorem \ref{main}
  says in particular that $\lambda_A=1$ if and only
  if $\lambda_B=1$.  Thus Theorem \ref{main}
  implies Theorem \ref{LR}.
  Now we prove Theorem \ref{main}.

  Let $f(A,B)=\lambda_A \lambda_B$.
    There exists a good pairs $(A_0,B)_0$ such that
    $f(A,B)=1$.  Take $A_0$ to be the regular $N$-sunburst and
  $B_0$ to be the suitably rotated copy.
  Call two good pairs $(A,B)$ and $(A',B')$ {\it closely related\/}
  if they differ only in the placement of a single line.
  For instance, we might have $A=A'$ and $B_k=B_k'$ except
  when $k=2$.   Any two good pairs can be connected by
  a finite sequence of closely related good pairs.
  In other words, we can get from one pair to the other by
  moving one line at a time.    
  In particular, we can start
  with $(A_0,B_0)$ and then reach an arbitrary good pair
  through a finite sequence of closely related pairs.
  For this reason, the next result implies Theorem \ref{main}.

  \begin{lemma}
    $f(A,B)=f(A',B')$ if $(A,B)$ and $(A',B')$ are
    closely related.
  \end{lemma}

  \startproof
  Given the invariance properties of $\lambda_A$ and $\lambda_B$
  discussed above, it suffices to prove our result in the special
  case already mentioned: $A'=A$ and $B_k=B_k'$ except when $k=2$.
  We can identify the space of complex lines through the origin
  with the Riemann sphere $\C \cup \infty$ in the usual way.
  We are just talking about the complex projective line here.
  Given $\zeta \in \C \cup \infty$ let
  $B(\zeta)$ be the $B$-sunburst obtained by replacing $B_2$ with
  the complex line $B_2(\zeta)$
  corresponding to $\zeta$.  Then there are
  parameters $\zeta,\zeta'$ such that
  $B=B(\zeta)$ and $B'=B(\zeta')$. We define
  $f(\zeta)=f(A,B(\zeta))$.
  
  There are $2$ bad values of $\zeta$ where $f$ is undefined, namely when
  $B_2(\zeta)=A_1$ or $B_2(\zeta)=A_3$.
  Given the nature of the construction,
  $f$ is a holomorphic function of $\zeta$, defined away from
  $2$ points on the Riemann sphere.

  Let us analyze the behavior of $f$ at the two bad values.
  We use the notation $g \sim h$ to denote the statement
  that the ratio $|g/h|$ is uniformly bounded away from
  both $0$ and $\infty$.   Here $g$ and $h$ are functions
which depend on the varying choice of $\zeta$.
  
  Suppose  that $B_2(\zeta)$ makes an angle of $\epsilon$ with $A_3$ and
  $|z_1|=1$.  Then $|z_3| \sim 1/\epsilon$.  The rest of the points
  are not affected much by the change.  This gives
  $\lambda_A \sim 1/\epsilon$.    A similar analysis shows
  that $\lambda_B \sim \epsilon$.  Hence
  $f$ is bounded in a neighborhood of
  the parameter $\zeta$ where
  $B_2(\zeta)=A_3$.  A similar analysis works for the other
  bad parameter.

  Since $f$ is a meromorphic function on the Riemann sphere
  with no poles, $f$ is constant.
  \endproof

  \subsection{Weaves}

  Let $(A,B)$ be a pair
  of $N$-sunbursts.
  For an even index $k$, we say that $(A,B)$ is {\it woven at\/} $k$
  if $B_{k}$ is parallel to a vector that starts at some interior point of
$A_{k-1}$ and ends at some interior point of $A_{k+1}$.
We call $(A,B)$ a {\it weave\/} if it is woven at $k$ for all
$k=2,4,...,2N$.
The pair $(A,B)$ in Figure 3.1 is a weave.
Our definition is more symmetric than it looks. 
Let $-B$ denote the sunburst obtained by
reflecting $B$ through the origin.

\begin{lemma}
  \label{switch}
  $(-B,A)$ is a weave if and only if $(A,B)$ is a
   weave.
\end{lemma}

\startproof
If $(A,B)$ is a weave, then
so is $(-A,-B)$.  We are just turning the picture upside down.
So, it suffices to prove the ``if'' direction.
Note that now the
weave property for $(-B,A)$ involves odd indices.

Let $\alpha=A$ and $\beta=-B$. Thus,
our new pair is $(\beta,\alpha)$.
The rays of $\beta$ are oriented
outward and the rays of $\alpha$
are oriented inward.  In particular
$\beta_k=B_k$ and $\alpha_k=-A_k$
for all relevant indices.

\begin{center}
\resizebox{!}{1.4in}{\includegraphics{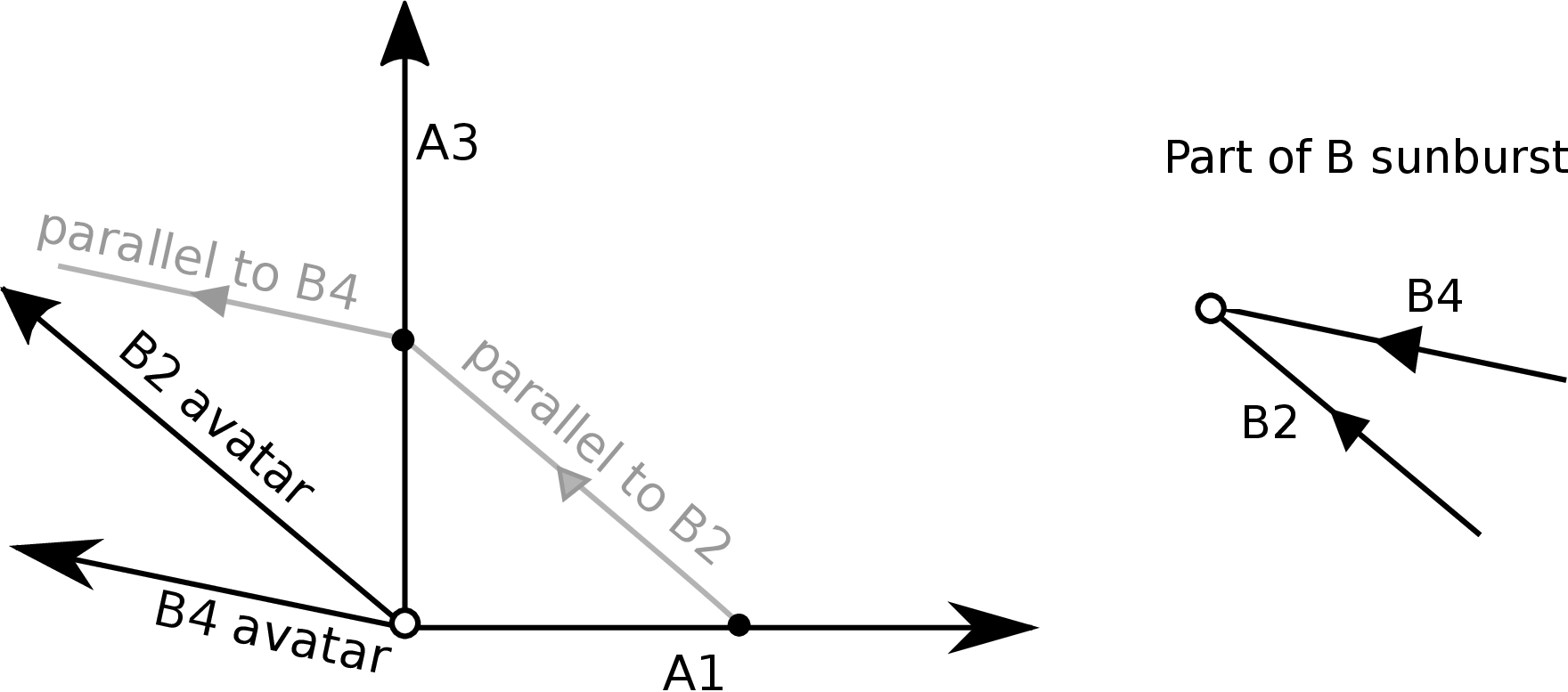}}
\newline
{\bf Figure 3.2:\/}  The relevant points and lines
\end{center}

We show that $(\beta,\alpha)$ is
woven at $3$.   The same argument works
for the other odd indices.
Say that an {\it avatar\/} of a ray is
a vector based at the origin and parallel to
the ray.
We normalize by an affine
transformation so that $A_1$ is the positive
$X$-axis and $A_3$ is the positive $Y$-axis.
Then, since $(A,B)$ a weave, $B_2$ has an avatar that
points into the $(-,+)$ quadrant.
Since $(A,B)$ is a weave, $B_4$ has an avatar
that points into the left halfplane.
Since $B$ is a sunburst with inwardly oriented
rays that progress counter-clockwise around as
the index increases, any avatar of $B_4$ lies
beneath any avatar of $B_2$.

The avatars of $\beta_2$ and $\beta_4$ equal
the avatars of $B_4$ and $B_4$ and
$\alpha_3$ is the negative $X$-axis.  From
what we have said above,
$\alpha_3$ is parallel to a vector connecting
a point on $\beta_2$ to a point on $\beta_4$.
Hence $(\beta,\alpha)$ is woven at $3$.
\endproof

 We call the system $(A,B')$ a {\it phase modification\/} of $(A,B)$
  if $B'$ is obtained by rotating $B$ about the origin by some angle.

\begin{theorem}
  \label{one}
  Let $(A,B)$ be a weave.  Then there is
  a  unique phase modification $(A,B')$ of $(A,B)$ which is a
  weave with periodic woven orbits.
\end{theorem}

The rest of this section is devoted to proving Theorem \ref{one}. 

\begin{lemma}
  \label{swap}
  Suppose $(A,B)$ is a weave.
  Then the orbit of $(a_1,b_2)$ is  woven.
\end{lemma}

\startproof
Since $(A,B)$ is a weave, the ray emanating
from $a_1$ and pointing in the direction of $B_2$ intersects
$A_3$.  Thus $a_3$ is well-defined.   According to our definition
in terms of particle, the transverse vector at $a_3$ points
into the sector bounded by $A_3$ and $A_5$.
Because $(B,-A)$ is also a weave, the ray
through $b_2$ and parallel to $-A_3$ intersects
$B_4$. Thus $b_4$ is well-defined.    According to our definition
in terms of particle, the transverse vector at $b_4$ points
into the sector bounded by $B_4$ and $B_6$.  Continuing
like this, we see that the forward orbit is woven.
\endproof

For each index, there is an open interval of phase modifications which
keep that part of the oriented weave condition.
The intersection $I$ parametrizes the phase modifications that are
weaves.

Suppose that $(A,B)$ is a weave.
Let $\lambda_A(A,B)$ be the function from
\S \ref{leftright}.   Since the orbits are woven, we have
$\lambda_A(A,B) \in (0,\infty)$. These orbits are
periodic iff
$\lambda_A(A,B)=1$.  We identify $I$ with $(0,1)$
and we orient $I$ so that as $t$ increases in $(0,1)$ the
corresponding sunburst $B'=B_t$ rotates clockwise.
We define $\lambda(t)=\lambda_A(A,B_t)$.
We want to see that there is a unique value $t \in (0,1)$ such
that $\lambda(t)=1$.  Theorem \ref{two} follows immediately
from the combination of the next two results.

\begin{lemma}
  \label{monotone}
  $\lambda$ is strictly increasing on $(0,1)$.
\end{lemma}

\startproof
Consider the orbit
$a_1(t),b_2(t),a_3(t),...$ from
Lemma \ref{swap}.  Here our notation
reflects the fact that this orbit dependn $t \in (0,1)$.
Notice that the ratio $\lambda_k(t)=\|a_{k+2}(t)\|/\|a_k(t)\|$
depends only on the triple of rays
$A_k, B_{k+1}(t),A_{k+1}$.  As $t$
increases, the ray $B_{k+1}(t)$ rotates
clockwise.  But then $\lambda_k(t)$ is
strictly  increasing.  Since
\begin{equation}
  \label{prod1}
  \lambda(t)=\lambda_1(t) \times ... \times \lambda_{2k-1}(t),
  \end{equation}
$\lambda(t)$ is also strictly increasing.
\endproof

\begin{lemma}
  $\lambda(t) \to 0$ as $t \to 0$ and
  $\lambda(t) \to \infty$ as $t \to 1$.
\end{lemma}

\startproof
We continue with the notation from the
proof of Lemma \ref{monotone}.
As $\lambda \to 0$, one of the $B$-rays
converges to one of the $A$-rays.
But then for the corresponding
index $k$, the quantity $\lambda_i(t)$ exits
every compact subset of $(0,\infty)$.
Since this quantity is decreasing, we see that
$\lambda_i(t) \to 0$.  At the same time,
the remaining factors in Equation \ref{prod1} do
not increase. Hence $\lambda(t) \to 0$.
The same kind of argument shows that
$\lambda(t) \to \infty$ as $t \to 1$.
\endproof

\subsection{A Calculus Interlude}

In this section we prove  an inequality that we use
in the next section.

\begin{lemma}
  \label{calc2}
  \label{calc3}
  Let $N \geq 4$ and
  $k=4,6,...,N$ and
  $y^*=\sin(\pi(k-2)/(2n))$.
    Then
  \begin{equation}
    \label{messy}
    ((k/2)-1) - y^*(N-(k/2)+1) <0.
    \end{equation}
\end{lemma}

\startproof
The function
$f(t)=\sin(\pi t)-t/(1-t)$ is positive on
$I=(0,1/2)$ because
$f(0)=f(1/2)=0$ and
$$f''(t)=-\pi^2 \sin(t) -\frac{2}{(1-t)^2} - \frac{2t}{(1-t)^3}<0.$$
Now,
let $t = (k-2)/2N$.  Note that
the conditions on $k$ give
$t \in (0,1/2)$.  The positivity of $f$ gives
$$y^* = \sin(\pi t) > \frac{t}{1-t}=\frac{k/2-1}{N-(k/2)+1}.$$
The last equality requires a bit of algebra.  Equation
\ref{messy} is a rearrangement of this inequality.
\endproof

\begin{center}
\resizebox{!}{1.1in}{\includegraphics{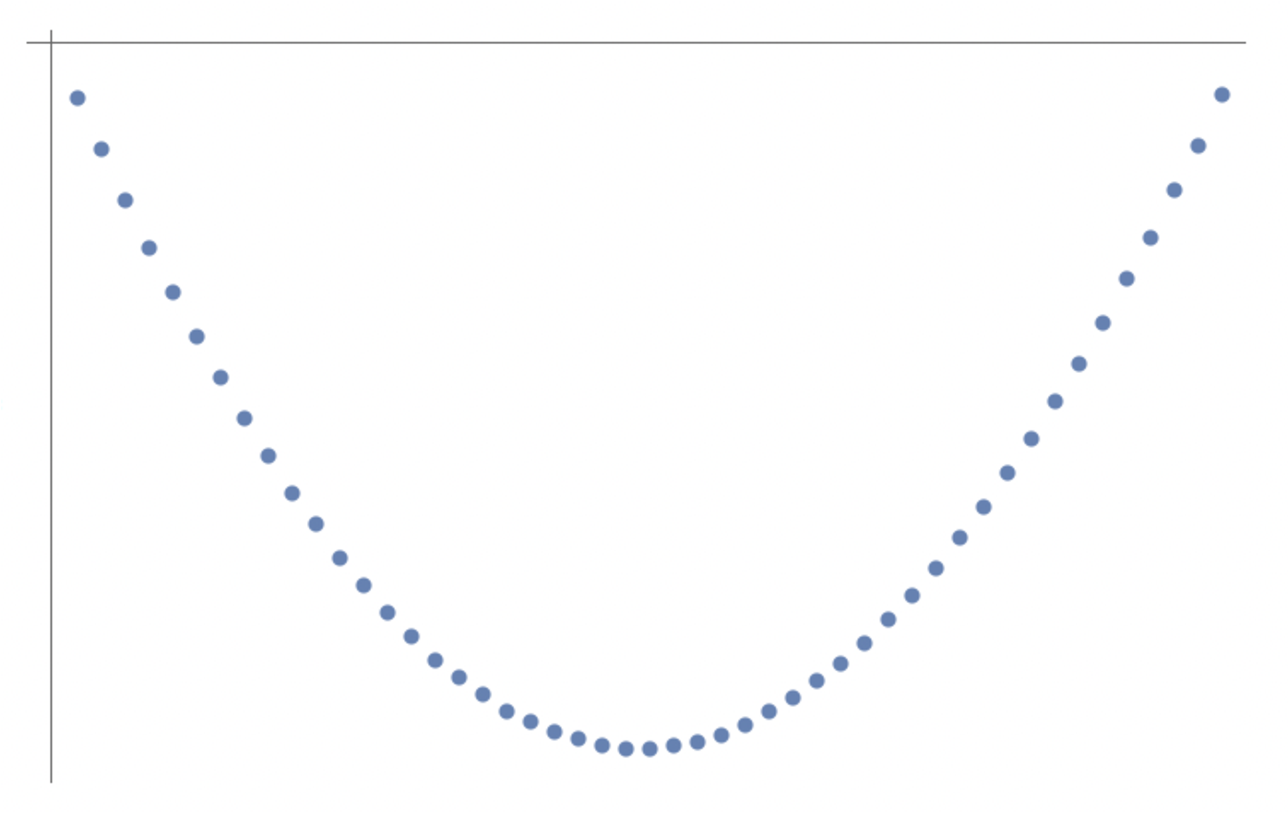}}
\newline
{\bf Figure 3.3:\/}  A plot of equation \ref{messy} for $N=100$ and $k=4,...,50$.
\end{center}

\subsection{Existence of Weaves}
\label{weaveexist}

Now we prove our restatement of
Theorem \ref{local}.

\begin{theorem}
  \label{two}
  If $A$ is regular and $B$ is balanced, then
  $(A,B)$ has a unique phase modification which is both a
  weave and has woven periodic orbits.
\end{theorem}

By Theorem \ref{one}, it suffices to prove that
$(A,B)$ has a phase modification which is a weave.
Let $S^1$ be the set of all
phase modifications $(A,B')$ of $(A,B)$.  For each even index $k$
there is an interval $I_k \subset S^1$ which parametrizes the
phase modofications that are woven at $k$.  We just need to prove that
$I=\bigcap I_k$ is nonempty.
Define
\begin{equation}
  \mu_N=\pi - \frac{2 \pi}{N}.
\end{equation}

\begin{lemma}
  Each interval $I_k$ has angular length $\mu_N$.
\end{lemma}

\startproof
Choose any point
$p \in A_{2k-1}$.  Then as $q \in A_{2k+1}$
moves from the origin to $\infty$, the vector
$\overrightarrow{pq}$ sweeps out an angle of $\mu_N$.
\endproof

When $B=A$, the
intervals $I_1,...,I_{2k-1}$ all coincide, by symmetry.
If $B$ is a general balanced $N$-sunburst,
we can find a homotopy $t \to B(t)$, through
balanced $N$-sunbursts, such that
$B(0)=A$ and $B(1)=B$.
For all indices $i,j$ we will show that the
relative displacement of $I_j(t)$ with respect to
$I_i(t)$ is less than $\mu_N$.
This means that all pairs of intervals intersect for all $t$.

\begin{lemma}
  Suppose that all the intervals $I_i(t)$ and $I_j(t)$
  intersect for all $t \in [0,1]$.  Then
  all the intervals $\{I_k(1)\}$ have a common intersection.
\end{lemma}

\startproof 
Recall a case of Helly's Theorem:
If a finite collection of open intervals in $\R$ pairwise
intersect, then their intersection is nonempty.  We think
of $\R$ as the universal cover of $S^1$.
We can lift our intervals to $\R$, so that for
$t=0$ they are all the same interval and the
lifts vary continuously with $t$.  But then the result about
relative displacement still holds, and all pairs of lifted
intervals intersect for all $t \in [0,1]$.  By Helly's Theorem,
all the lifted intervals intersect (in particular) for $t=1$.
Pushing the intersection point
down to $S^1$, we see that all the intervals
$\{I_k(1)\}$ also intersect.
\endproof

To finish our proof of Theorem \ref{two} we need to establish
our claim about the relative displacement of pairs of intervals.
Without loss of generality, it suffices to consider the
intervals $I_0(t)$ and $I_k(t)$ for $k=2,4,...,2N-2$.
The relative displacement of our intervals does not change if
we post-compose our homotopy $B(t)$ with an arbitrary
continuous family of rotations. So, it suffices to
consider the case when $A_1$ is the positive $X$-axis
and $B_0(t)$ is the positive $X$-axis for all $t \in [0,1]$.

Let $\theta_k(t)$ denote the angle between $B_k(t)$ and
the positive $X$-axis.  We have $\theta_0(t)=0$ for all $t$
and $\theta_k(0)=\pi k/N$ for all $k=0,2,4,...,2N-2$.
The relative displacement between $I_k(1)$ and $I_k(0)$ is
$|\theta_k(1)-\theta_k(0)|$.  We just have to prove that
this quantity is less than $\mu_N$ for all even $k$.

  \begin{lemma}
    $\theta_k(1)-\theta_k(0)<\mu_N$ for all $k=2,...,2N-1$.
  \end{lemma}

  \startproof
  We will assume that there is an index $k$ such that
  $\theta_k(1)-\theta_k(0) \geq \mu_N$ and
  derive a contradiction.  By assumption
  $$\theta_k(1) \geq \theta_k(0) + \pi - \frac{2\pi}{N} =
  \pi + \frac{k \pi}{N} -\frac{2 \pi}{N} =
  \pi +\frac{(k-2)\pi}{N}.$$

  We first treat two special cases.
  \begin{itemize}
  \item When $k=2$ we have $\theta_2(1) \geq \pi$.
    This contradicts the fact that $\theta_2(1)<\pi$ is the angle
    between $B_0(1)$ and $B_2(1)$.
\item   When $k \geq N+2$ we have
  $\theta_k(1) \geq 2\pi$.  This contradicts the fact that
$\theta_k(1)<2\pi$ for $k=2,4,...,2N-2$.
\end{itemize}

So, we only need consider the cases
$k=4,6,...,N$.  Compare Lemma \ref{calc2}.
  Let $\beta_j$ denote the unit vector along
  $B_j(1)$.  We rotate the picture about
  the origin so that $\beta_0$ and $\beta_k$ have
  the same negative $y$-coordinate, namely
  $-y^*$, the quantity from Lemma \ref{calc2},
  with $\beta_0$ being on the right and
  $\beta_k$ being on the left.
  The $(k/2)-1$ unit vectors $\beta_2,...,\beta_{k-2}$
  all have $y$-coordinates less or equal to $+1$.  The remaining
  $N-(k/2)+1$ unit vectors,
  being between $\beta_k$ and $\beta_0$ in the
  circular order, all have $y$ coordinates less or
  equal to $-y^*$.   The average of the $y$-coordinates
  of all $N$ unit vectors is at therefore at most
  $$(k/2-1) - y^* (N-(k/2)+1)<0.$$
  The inequality comes from Equation \ref{messy}.
  This contradicts the fact that $B(1)$ is balanced.
  \endproof

  \begin{lemma}
    $\theta_k(0)-\theta_k(1)<\mu_N$ for all $k=2,...,2N-1$.
  \end{lemma}

  \startproof
  This follows from the previous result and reflection symmetry.
  Let $\overline B$ denote the balanced $N$-sunburst obtained by
  reflecting $B$ in the $X$-axis and dihedrally relabeling the rays
  so that they again go counterclockwise around the origin.  Applying
  the previous result to $\overline B$ and we get
    $\overline \theta_k(1)-\overline \theta_k(0)<\mu$ for all
  $k=2,..,2N-2$.  The homotopy we take is just the original homotopy
  conjugated by reflection in the origin, with the points dihedrally relabeled.
  At the same time, we have
  $\theta_k(t)=2\pi - \overline \theta_{2N-k}(t).$
  Hence
    $$\theta_k(0)-\theta_k(1) =
  -\overline \theta_{2N-k}(0)+\overline \theta_{2N-k}(1)=
  \overline \theta_{2N-k}(1)-\overline \theta_{2N-k}(0)<\mu.$$
This completes the proof.
  \endproof

  \noindent
  {\bf Remark:\/}
  One might wonder about relaxing the hypotheses of
  Theorem \ref{two}.   For instance, is the result true
  for $(A,B)$ if both $A$ and $B$ are regular?    In the
  even case, at least, one needs extra hypotheses.
  When $N$ is even and $A$ is not centrally symmetric,
  there always exists a centrally symmetric (and hence balanced)
  sunburst $B$ such that $(A,B)$ has no phase
  modification which is a weave.  The sunburst $B$
  consists of two clusters of $N/2$ nearly identical lines
  which are diametrically opposed from each other.
  We leave the details of this to the interested reader.
  In the odd case I have not been able to think of an
  easy counter-example like this.

\newpage

\section{Correspondence between Polygons}

\subsection{The Main Construction}

We consider planar equilateral polygons in which all
the sides have the same length.  We also consider
planar equiangular polygons in which all the interior
angles are the same.  We consider these polygons
up to the equivalence of orientation preserving
similarity.  In both cases we only consider strictly
convex polygons which are oriented counterclockwise
around the regions they bound.  In the equilateral
case we normalize so that the sides have
unit length.  In the equiangular case we 
normalize so that the sides are parallel to the
relevant roots of unity.
\newline
\newline
\noindent
{\bf Main Construction:\/}
Let $L$ be a strictly convex $N$-gon with unit vector edges
$\beta_2,\beta_4,...,\beta_{2N}$.  Since $L$ is closed, we have
$\beta_2+ ... + \beta_{2N}=0$.
Let $B_k$ be the ray through the origin that starts at $0$ and
contains $\beta_k$.
The fact that $L$ is strictly convex means that
$B_2,...,B_{2N}$ is an $N$-sunburst.
By construction, $B$ is a balanced sunburst.
A different representative of $[L]$ would give rise
to a balanced sunburst $B'$ which is a rotation of $B$.
Let $A$ be the regular $N$-sunburst.
Applying Theorem \ref{two}, we know that
there is a unique phase modification
$(A,B')$ which is a weave and which has periodic orbits.
We choose any periodic orbit and consider the
convex $N$-gon $P_L$ which is the image that
the orbit traces around $B$.   Our association
is \begin{equation}
  \label{assoc}
  [L] \to [P_L].
\end{equation}
The labeling is such that the $k$th vertex of $P_L$ lies
on the edge $B_k$.  
By construction, this association
is well-defined, independent of choices.
\newline

\noindent
{\bf Remarks:\/}
(1)   
The construction above only works in the
strictly convex case, but we discuss the general
case in \S \ref{all}.
\newline
(2)
Spherical duality interchanges equiangular and equilateral
spherical polygons.   One might wonder if one could get a
different correspondence like Equation \ref{assoc} by taking
a limit of this process as the sphere radius tends to infinity.
(One of the referees posed this question.)
I think that this probably will not work.
Consider what happens when we have
a unit-sized regular polygon on a growing sphere. Then
the diameter of the dual polygon tends to infinity and
it seems impossible to extract a Euclidean limit.
The same problem would happen for other polygons.

\subsection{Computing the Correspondence}
\label{compute}

In the introduction we claimed that the correspondence in
Equation \ref{assoc} is easy to compute.   Figures 1.3 and
3.1 come from our computer program which does the
computations.  Here we discuss the method.
Once we have the pair $(A,B)$ we find the intervals.
$I_1,...,I_N$ where $I_k \subset S^2$ is the interval
of rotations $R_{\theta}$ such that
$(A,R_{\theta}(B))$ is $k$-intervoven. See \S \ref{weaveexist}.

Next, we compute $I=\bigcap I_k$.   We have the
map $h: I \to \R_+$ which
computes the holonomy $h(\theta)$
of the orbits associated
to the pair $(A,R_{\theta}(B))$ for $\theta \in I$.
We want to solve the equation $h(\theta)=1$.
We use the bisection method.  We choose
two values $\theta_1,\theta_2$ respectively
very near the two boundary points of $I$.
We then perform the following bisection algorithm.
\begin{enumerate}
  \item Start with a pair $(\theta_1,\theta_2)$ such that
    $h(\theta_1)<1<h(\theta_2)$.
  \item Let $\theta_3=(\theta_1+\theta_2)/2$.
  \item If $h(\theta_1)>1$ replace
    $(\theta_1,\theta_2)$ with
    $(\theta_1,\theta_3)$.
  \item If $h(\theta_1) < 1$ replace
            $(\theta_1,\theta_2)$ with
    $(\theta_3,\theta_2)$.
  \item Return to Step 1 using the smaller angle  interval.
    \end{enumerate}
      If we iterate this $M$ times we get the
      correct value of $\theta$ up to about
      $2^{-M}$.
      Once we have our good approximation of
      $\theta$ it is a simple matter to actually
      compute the symplectic tiling orbit.
      \newline
      \newline
      {\bf Remark:\/}
The construction in [{\bf KM\/}] starts with
an equiangular $N$-gon, and then take the Riemann map
to the unit disk,  This gives them $N$ unit complex numbers
$u_1,...,u_N$.
They then use the fact (Springborn's Theorem [{\bf Sp\/}]) that there is
a unique point in the hyperbolic plane such that a
Moebius transformaton $M$ mapping this point to the origin
carries $u_1,...,u_n$ to unit complex numbers whose sum is $0$.
(This part of the construction is similar to ours.)
The numbers $u_1',...,u_n'$ are then interpreted as the
direction vectors of a unit equilateral $N$-gon.
Here we have set $u_k'=M(u_k)$.  This construction
produces a unique equilateral $N$-gon up to scale.
This method requires
      something like the computation of the Riemann map.
      There are several methods for doing this for the
      case at hand.  Concretely, one could numerically
      integrate the Christoffel transform.   Alternatively,
      one could use elegant circle-packing methods
      [{\bf St\/}].  All these methods seem very
            computationally involved.
      
\subsection{Bijective Nature of the Correspondence}

We fix some integer $N \geq 3$ and consider
all our constructions relative to $N$.  In this section
we prove that our association in Equation \ref{assoc} is
a bijection.   We explain first how to construct the
inverse map and then we justify the assumptions needed
for the construction to work.

Let $[P]$ be an equivalence class of convex equiangular $N$-gons.
Let $P$ be some representative.  Let
$P_2,P_4,...,P_{2n}$ be the vertices of $P$.
Say that a point $p$ in the region
bounded by $P$ is a {\it balance point\/} if the sum
$\sum \beta_k=0$, where $\beta_k$ is the unit vector parallel
to $P_k-p$.    We will show below that there is a unique
balance point.  This result is similar in spirit to
Springborn's Theorem [{\bf Sp\/}].

Given our balance point $p$, we
let $B$ be the balanced sunburst defined by
the vectors $\beta_2,...,\beta_{2n}$ we have just defined.
Finally, let $L$ be the equilateral $N$-gon whose successive edges are
parallel to these vectors.  The balance condition guarantees
that $L$ is closed and the sunburst condition guarantees
that $L$ is strictly convex.
Our construction shows how to recover $[L]$ from $[P_L]$, and
this shows that our association is injective.
At the same time, our construction shows that our association
is surjective.  Hence, our association in
Equation \ref{assoc} is bijective.

\begin{lemma}
  A strictly convex $N$-gon has at most one balance point.
\end{lemma}

\startproof
For this proof we do not use the equiangular property.
We will suppose this is false and derive a contradiction.
Let $P$ be a strictly convex polygon which supposedly
has at least $2$ balance points.  We rotate the picture so
that both balance points $p_1,p_2$ lie on the $X$-axis
and $p_1$ is on the left.

\begin{center}
\resizebox{!}{1.7in}{\includegraphics{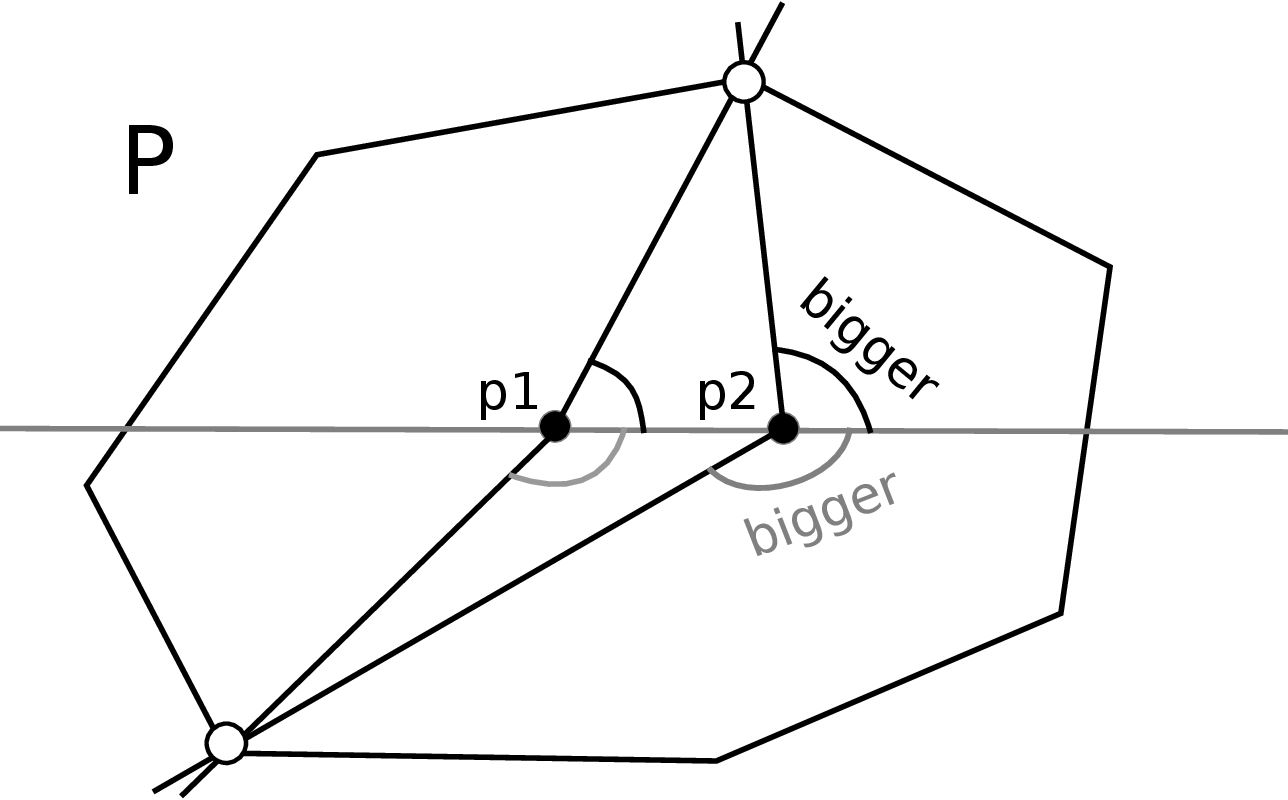}}
\newline
{\bf Figure 4.1:\/}  The angles with the $X$-axis.
\end{center}

For each $i=1,2$ let $\{\beta_{ij}\}$ denote the set of
  unit complex numbers that are parallel to the
  vectors $P_j-p_i$.  The basic property is that
  the angle that the vector $\beta_{1j}$ makes with the
  $x$-axis is less than the angle that $\beta_{2j}$ makes
  with the $x$-axis.   Figure 4.1 shows this in action.
  From this we see that the center of mass of
  $\{\beta_{2j}\}$ must lie to the left of the center
    of mass of $\{\beta_{1j}\}$, contradicting the claim
    that both centers of mass are the origin.
    \endproof

    \begin{lemma}
      \label{exist}
  A strictly convex equiangular $N$-gon has a balance point.
\end{lemma}

\startproof
This proof  uses the equiangular property.
When $N=3$ the polygon must be an equilateral
triangle, and then the center of symmetry does the job.
Likewise, when $N=4$ the polygon must be a rectangle,
and again the center of symmetry does the job.  So,
we take $N \geq 5$.

Let $P$ be a strictly convex equiangular $N$-gon.   We consider a
simply connected domain $D$ in the plane as follows.
We start with the closure of the domain bounded by $P$ and then
we chop off small isosceles triangular
neighborhoods of the vertices of $P$.
Figure 4.2 shows the picture.
The boundary of $D$ is a convex $2N$-gon which equals
$P$ at most places but then takes small ``shortcuts'' into
the interior of the region bounded by $P$ near the vertices.

For each $p \in D$ we consider the vector
$V_p = \sum \beta_{p,j}$, where $\beta_{p.j}$ is
the unit vector parallel to the one pointing from $p$ to $P_j$.
We are looking for a place where $V_p$ vanishes.
If suffices to prove that $V_p$ points inward at
$\partial D$.  (As the proof develops, we will explain
more precisely what this means.)
The boundary $\partial D$ has two kinds of points, those which also
lie in $P$ and those which do not.
We consider these kinds of points in turn.

\begin{center}
\resizebox{!}{1.4in}{\includegraphics{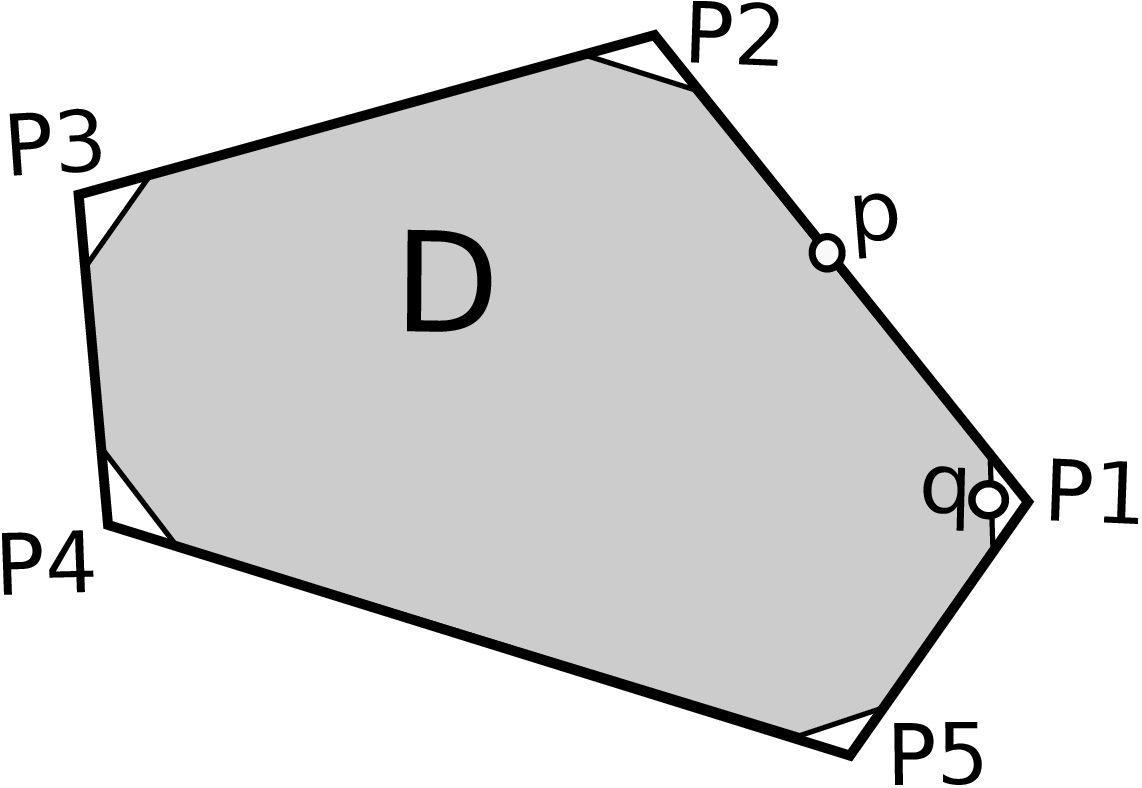}}
\newline
{\bf Figure 4.2:\/}  The domain $D$ and two kinds of boundary points
\end{center}

The point $p$ in Figure 4.2 also lies in $P$.
For points like $p$, all but two of the
vectors $\beta_{p.j}$
point into $D$ and the other two point along $\partial D$.
Let us say this a bit more formally.
The $p$ lies in an edge $e_p$ of $P$,
and $e_p$ is contained in a line $L_p$. The line $L_p$ bounds a
halfplane $H_p$ that contains the other vertices of $P$, and
$\beta_{p,j}$ points into $H_p$ for all but $2$ vertices.
For the other two vertices, namely the vertices of $e_p$,
the two unit vectors in question point along $L_p$.
Adding up all these unit vectors, we see that $V_p$
points into $H_p$.  This is to say that $V_p$ points
into $D$.

The point $q$ in Figure 4.2 does not lie in $P$.
We rotate the picture so that the edge containing
$q$ is vertical, as shown in Figure 4.2.
We also relabel so that $q$ is near $P_1$.
The $x$-coordinate of $\beta_{p,1}$ is at most $1$.
The remaining vertices lie on lines which make an
angle of at least $\pi/N$ with the vertical line
through $P_1$.  Here we are using the equiangularity
condition.     But this means that the sum of the
$x$-coordinates of our remaining vectors is
at most
\begin{equation}
  \label{negsum}
  -(N-1) \times (\sin(\pi/N) - \epsilon).
  \end{equation}
Here $\epsilon$ is a number can make as small as
we like by controlling the size of the isosceles
neighborhoods we used in defining $D$.
As long as we take $\epsilon$ sufficiently small,
the number in Equation \ref{negsum} is
always less than (meaning more negative than)
$-2$.    Hence the vector $V_q$ has negative
$x$-coordinate.  This shows that $V_q$ points
into $D$.

Now we know that our vector field $p \to V_p$ is inward-pointing
on $\partial D$.  A well known result about the index of vectorfields
now shows that $V_p$ vanishes somewhere in the interior of $D$.
Hence $P$ has a balance point.
\endproof

\subsection{Algebraic Nature of the Correspondence}
\label{algXX}

In this section we explain the sense in which the
association in Equation \ref{assoc} is algebraic.
The result here feeds into the proof of Theorem
\ref{algebra}.  

\begin{lemma}
  \label{alg}
  If $[L]$ is algebraic then $[P_L]$ is algebraic.
\end{lemma}

\startproof
When $L$ is algebraic, the rays $B_2,B_4,...,B_{2N}$
are also algebraic.  For instance, they have
algebraic slopes.  The holonomy $\lambda$
of an orbit
associated to $(A,B)$ is an algebraic function of
these slopes.    To understand the phase-modification
part of the construction we choose a rational
parametrization of the circle, as in
Equation \ref{RAT}.    For each $t \in \R \cup \infty$
the slopes of the $t$-rotated sunburst $B'$ are
rational functions in $t$, with algebraic coefficients.

We think of the holonomy $\lambda(t)$ as a function of
the rotation parameter $t$.  The function $\lambda(t)$ is
also a rational function with algebraic coefficients.
Setting $\lambda(t)=1$ and solving, we see that
the choice of $t$ which makes $(A,B')$ have periodic
orbits is algebraic.  But then if we scale so that one of the
vertices of $P_L$ is algebraic then all the vertices will
be algebraic.
\endproof

\begin{lemma}
  \label{algebra2}
  $[P_L]$ if and only if $[L]$ is algebraic.
  \end{lemma}

  \startproof
  In view of Lemma \ref{alg}, we just need to prove
  that $[L]$ is algebraic when $[P_L]$ is algebraic.
  Let $P=P_L$ be an algebraic representative.
  We claim that the balance point of $P$ is algebraic.
  (I am grateful to Joe Silverman for supplying the proof.)

We will use complex notation.
Let $P= (p_1,...,p_N) \in \C$. 
For each sequence
$\epsilon=(\epsilon_1,...,\epsilon_N) \in \{\pm 1\}^N$ define
\begin{equation}
  F_{\epsilon}(P,z)=\sum_{i=1}^N \epsilon_i \frac{z-p_i}{\overline z-\overline p_i} = 0.
\end{equation}
The balance point solves the equation
$F_{\epsilon}(P,z)$ when $\epsilon=(1,...,1)$.
The product
\begin{equation}
  G(z,p)=\prod_{\epsilon \in \{\pm 1\}^N} F_{\epsilon}(z,p)
\end{equation}
is unchanged if we change the signs of any subset of the
square roots in the last equation.   Hence $G(P,z) \in \Q(P,z)$,
the ring of rational functions in $z$ and $P$.  Hence
the roots of $G(P,z)$, for algebraic $P$, are also algebraic.
The balance point is one such root.

Since the balance point is algebraic,
the rays describing the $B$ sunburst are
algebraic.  So, if we start with $a_1$ and $a_2$ algebraic,
the whole orbit remains algebraic.
\endproof

\noindent
{\bf Remark:\/}
The algebraic structure of the balance point might be
quite complicated.  Consider the modest example of an
integer pentagon with vertices
$$(0,0), \hskip 15 pt
(1,0), \hskip 15 pt (2,2), \hskip 15 pt
(1,2) \hskip 15 pt (0,1).$$
Peter Doyle played around with this
in Mathematica and found that the minimal polynomial
for the first coordinate of the balance point has
degree $48$ and the coefficients mostly have about
$30$ digits.

\newpage

\section{Hyperbolic Structure}

\subsection{The Thurston Construction}
\label{hypXX}

William Thurston's paper [{\bf T\/}] constructs
complex hyperbolic structures on spaces of
flat cone spheres. See my notes [{\bf Sch1\/}]
for an exposition of [{\bf T\/}].
A special case of a flat
cone sphere is the double of a convex equiangular
polygon.   The corresponding subspace sits as a
totally real slice of the convex hyperbolic moduli
space.  This imparts a real hyperbolic structure
on the space of convex equiangular $N$-gons.

In this section I will give an elementary account of
the construction which does not go through complex
hyperbolic geometry. It is possible that I learned this
construction from Thurston when I was a graduate
student at Princeton University and it is also possible
that I worked it out myself sometime later.
There are a number of similar
accounts in the literature.  See e.g.
[{\bf BG\/}] for the general case and
[{\bf Cal\/}] for the pentagonal case.
\newline
\newline
{\bf Linear Coordinates:\/}
We start
with $N$ parallel familes of lines, with each family being
parallel to a different $N$th root of unity.  These families
are cyclically ordered, according to the roots of unity.
We interpret an equiangular $N$-gon as a selection
$\ell_1,...,\ell_n$ of lines, one from
each family.  The vertices of the $N$-gon
are given by $\ell_1 \cap \ell_2$, $\ell_2 \cap \ell_3$, etc.
This interpretration gives a natural identification
of the space of equiangular $N$-gons with $\R^N$.
To get a concrete coordinatization we could pick
some line $L$ in the plane, not parallel to any of
the families, and then use the intersection
$\ell_1 \cap L,...,\ell_N \cap L$ give $N$
linear coordinates.  In other words, we
are identifying $\R^N$ with $L \times ... \times L$.
A different choice of $L$ would give us a linear
change of coordinates.

We now mod out by translations.  This identifies
the space of equiangular $N$-gons mod isometry
with $\R^{N-2}$.   Figure 5.1 shows, in the pentagon
case, how we can introduce concrete coordinates
on $\R^{N-2}$ which are linear functions of
the coordinates discussed above.

\begin{center}
\resizebox{!}{1.2in}{\includegraphics{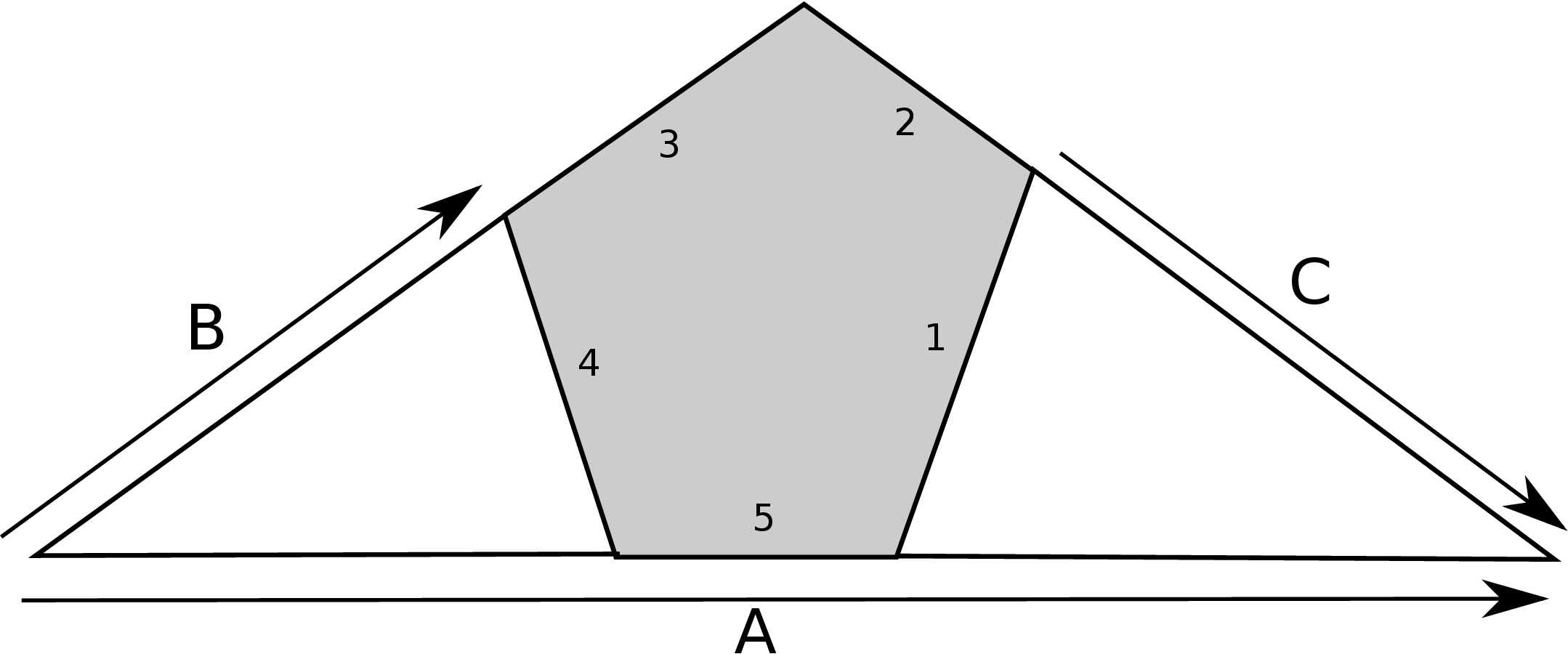}}
\newline
{\bf Figure 5.1:\/}  Coordinates on the space of equiangular pentagons
\end{center}

It is important to emphasize that these coordinates
$A,B,C$ are signed distances.  They look nice in the
convex case but they are defined even in the non-convex
cases.  For instance, in vector notation,
$$A=\big((\ell_2 \cap \ell_5) - (\ell_3 \cap \ell_5) \big)\cdot (1,0).$$
Also, it is important to note that any system of coordinates based
on a similar contruction (with other choices) would result in
a linear change of variables implemented by a matrix
with algebraic entries.

Figure 5.2 shows similar coordinates for the case of hexagons
and $7$-gons.  We had to make some choices to get these
coordinates, but any similar system
would be related by a change of coordinates implemented by
an algebraic matrix.

\begin{center}
\resizebox{!}{1.25in}{\includegraphics{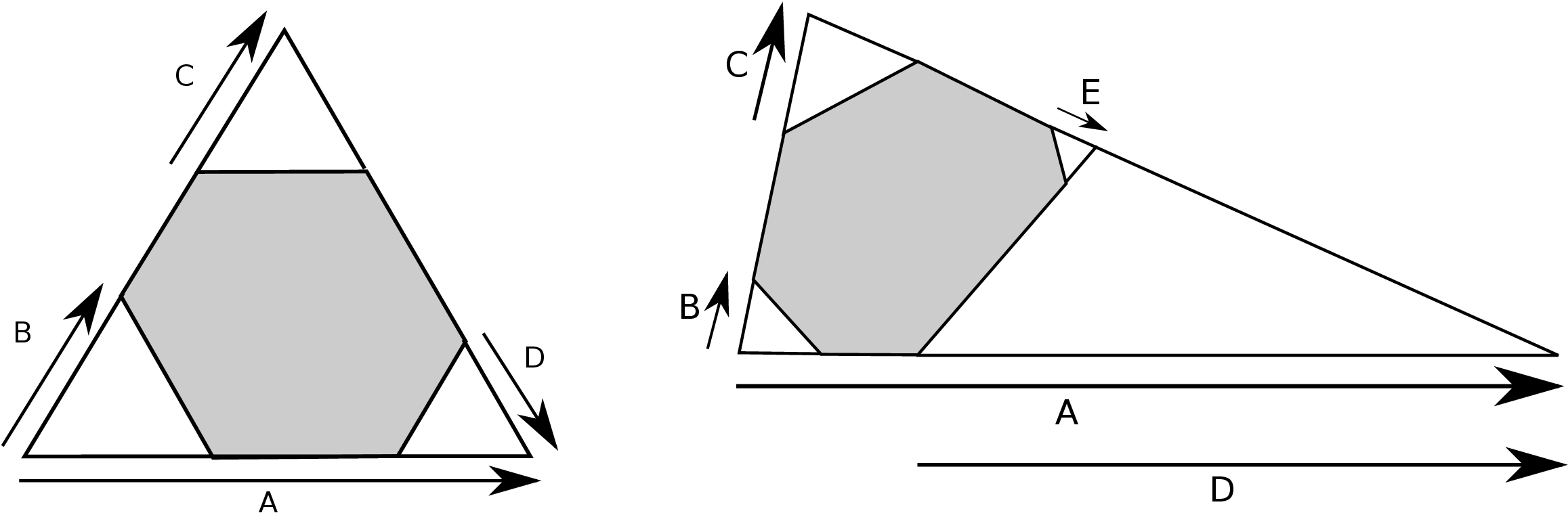}}
\newline
{\bf Figure 5.2:\/}  The case of hexagons and $7$-gons.
\end{center}

\noindent
{\bf The Signed Area:\/}
Now we consider the signed area in these coordinates.
In the pentagonal case we have
\begin{equation}
  \label{area}
  {\rm area\/} = -\alpha A^2 - \beta B^2  + \gamma C^2,
\end{equation}
where $\alpha,\beta,\gamma$ are positive constants that
do not depend on the choice of pentagon.
Geometrically, the area of the pentagon is the
area of the big triangle minus the area of the two small
triangles. The area of the big triangle is a quadratic function
of $C$ and the constant only depends on the shape of the
triangle.  Likewise the areas of the smaller two triangles are
quadratic functions in $A$ and $B$ with the same properties.

For hexagons and $7$-gons, Equation \ref{area} would
respectively have $4$ and $5$ quadratic terms with
constant coefficients and all but one being negative.
In the general case there would be $N$ quadratic terms with
constant coefficients with all but one being negative.
Speaking more abstractly,
the area of the $N$-gon given by
$V=(A,B,C,D,...)$ has the form
$Q(V,V)$ where $Q$ is a quadratic form of signature $(1,N-3)$.
\newline
\newline
{\bf The Lorentz Model:\/}
Now, we are interested in the space of {\it equivalence classes\/} of
equiangular $N$-gons.  Up to translation, we can get a unique
representative of each equivalence class by scaling so that the
area is $1$.  But then we can identify our space of equivalence
classes with one sheet of the hyperboloid in $\R^{1,N-3}$ given by
$Q(V,V)=1$.   This is a well-known {\it Lorentz model\/}
of $\H^{N-3}$, hyperbolic
space of dimension $N-3$.
\newline
\newline
{\bf Interaction with Convexity:\/}
Let ${\cal C\/}_N$ be the domain in
$\H^{N-3}$ corresponding to the
space of strictly convex equiangular $N$-gons.
In general, ${\cal C\/}_N$ is the interior of a convex polyhedral
domain.   The points on the boundary of
${\cal C\/}_N$ correspond to
degenerate polygons in which one or more
edge has collapsed to a point.

For instance, in the case of pentagons,
the boundary of ${\cal C\/}_5$ has
$5$ edges and $5$ vertices.
Referring to Figure 5.1, two of the edges
correspond to $B=0$ and $C=0$, and their
vertex intersection corresponds to $B=C=0$.
There is an elegant way to see the geometry
of ${\cal C\/}_5$.
Figure 5.3 shows the
{\it butterfly move\/} $B_2$ for pentagons.

\begin{center}
\resizebox{!}{3in}{\includegraphics{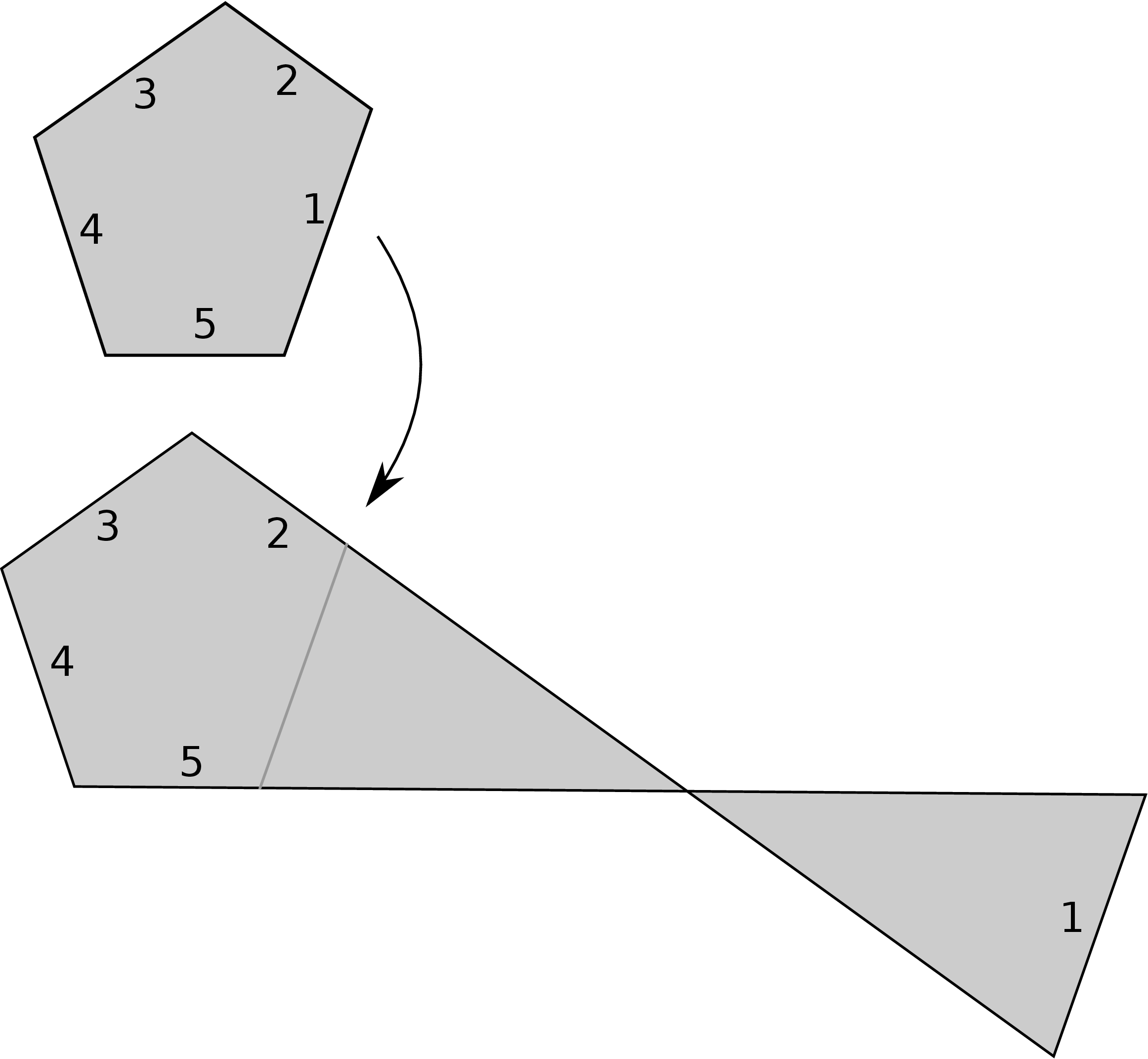}}
\newline
{\bf Figure 5.3:\/}  The butterfly move $B_2$.
\end{center}

The $5$ fixed lines in the hyperbolic
plane $f_1,f_3,f_5,f_2,f_4$ are consecutively
perpendicular in the cyclic sense.
The reason they are perpendicular is that
the corresponding butterfly moves commute! 
These fixed lines are the extensions
of the edges of a regular right-angled pentagon.
The interior of this pentagon is the space of
convex unit equiangular pentagons modulo
similarity.

The case of hexagons is also possible to understand.
In this case ${\cal C\/}_6$ has $6$ faces and $5$
vertices.  Two of the vertices lie in $\H^3$.  These
correspond to the two equilateral triangles we get
by collapsing either the even or the odd edges of
our hexagon.   The other three vertices are
ideal vertices.  These correspond to the hexagons one
gets by letting a pair of opposite sides get very long.
Figure 5.4 shows what we mean.

\begin{center}
\resizebox{!}{.9in}{\includegraphics{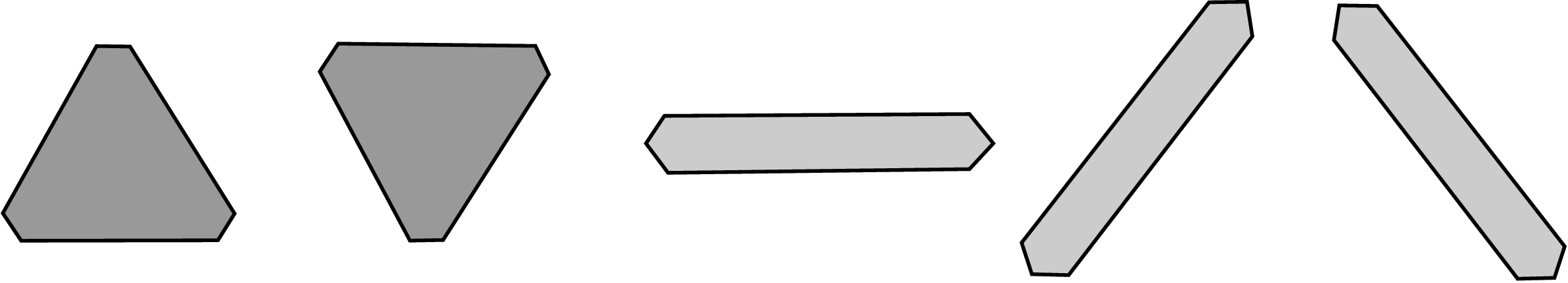}}
\newline
{\bf Figure 5.4:\/}  Hexagons near each of the $5$ vertices of ${\cal
  D\/}_6$.
\end{center}

Whenever two of these faces intersect, the corresponding
butterfly moves commute. Thus, all the faces which meet
do so at right angles.  The domain ${\cal C\/}_6$ in fact
is the interior of a triangular bi-pyramid.  Each half of the
triangular-bi-pyramid is obtained by coning one face of
a regular ideal octahedron to the center of mass.  This
half is a pyramid, one of whose faces corresponds to the
octahedron.  Call this the {\it blue face\/}.   Call the
other faces the {\it red faces\/}.  The red faces meet at
right angles and each red face meets the blue face at
an angle of $\pi/4$.  We get ${\cal C\/}_6$ by gluing
two of these pyramids together across their blue faces
and we are left with the $6$ red faces.

In general, if $P=\ell_1,..,\ell_N$ then
$B_k(P)=\ell_1',...,\ell_N'$ where
$\ell_j'=\ell_j$ for all $j \not = k$ and
the two lines $\ell_k, \ell_k'$ are equidistant
from the intersection $\ell_{k-1} \ell_{k+1}$.
The operation $B_k$ is linear and preserves
signed area.  Hence $B_k$ is a Lorenz transformation
and induces a hyperbolic isometry on the hyperbolic
structure we have explaned.  Moreover $B_k$ is an
involution.  The fixed point set of $B_k$ is the set of
all degenerate polygons in which $\ell_{k-1}, \ell_k,\ell_{k+1}$
have a common point.  This is a codimension one set.
Hence $B_k$ is a hyperbolic reflection.
Note that $B_{a}$ and $B_b$ commute as long as
$a,b$ are not cyclically consecutive.  The corresponding fixed
point sets are perpendicular hyperplanes.

\subsection{Putting it Together}
\label{put}

For each $N \geq 5$, the Thurston construction produces an open
polyhedral convex domain ${\cal C\/}_N \subset
\H^{N-3}$ whose interior
parametrizes the equivalence classes of
strictly convex equiangular $N$-gons.
Now we get to the punchline.
Using our correspondence
from Equation \ref{assoc}
we get the same hyperbolic structure on the space of
equivalence classes of strictly convex equilateral polygons.
\newline
\newline
\noindent
{\bf Proof of Theorem \ref{algebra}:\/}
Suppose that $L$ is an equilateral polygon with algebraic vertices.
Then by Lemma
\ref{algebra2} we can find a representative $P_L$
having algebraic vertices.   When we
scale $P_L$ to have unit area we are
scaling by an algebraic number, so we can
take $P_L$ to have unit area.
But then our special coordinates for $P_L$ are
also algebraic.  Any choice of special coordinates
would have this property, because they all differ
by the action of an algebraic linear matrix.
So, the coordinates of $P_L$ in Lorentz space
$\R^{1,N-3}$ are algebraic.  This is the same
as saying that the coordinates in $\H^{N-3}$ are algebraic.
Conversely, if $P_L$ has algebraic coordinates in
$\H^{N-3}$, then we can take a representative of
$P_L$ such that the special coordinates taken above
are all algebraic and one of the vertices is algebraic.
But then when we reconstruct $P_L$ from a single
vertex and from the coordinates we get an algebraic
polygon.  But then, by Lemma \ref{algebra2} again,
$[L]$ is also algebraic.
\endproof

\noindent
{\bf Remark:\/}
Finding 
the coordinates in $\R^{1,N-3}$ of a given class of
equiangular polygon is straightforward.
We just choose our coordinates as above and
compute. Thus, in view of the discussion in
\S \ref{compute}, it is easy to accurately estimate
the point in $\H^{N-3}$ which parametrizes a
given equilateral $N$-gon.

\subsection{Beyond the Convex Case}
\label{all}

Now we consider general equilateral polygons.
First of all, we widen our equivalence relation so
that two polygons are equivalent if and only if
there is a similarity which maps one to the other.
The similarity here need not be orientation preserving.
If we restrict our attention to the strictly convex case,
this widening of the equivalence relation changes
nothing, because the counterclockwise-oriented convex
polygons we have been considering above
are equivalent in the wider sense if and only if
they are equivalent in the narrow sense.

  The group $S_N$ of permutations
  acts naturally and continuously on the moduli space of equilateral $N$-gons:
We can encode an equilateral $N$-gon by an ordered list
$e_1,...,e_N$ of unit vectors.   Given a permutation
$\pi \in S_N$ we get the new list
$e_{\pi(1)},...,e_{\pi(N)}$ of edges, and we can build a
unique equivalence class of $N$-gon which corresponds
to this list.  The only thing we need from our vectors
is that they sum to zero, and this is unchanged by
permutation.  Figure 5.5 shows this in action for the
regular pentagon.

\begin{center}
\resizebox{!}{1.2in}{\includegraphics{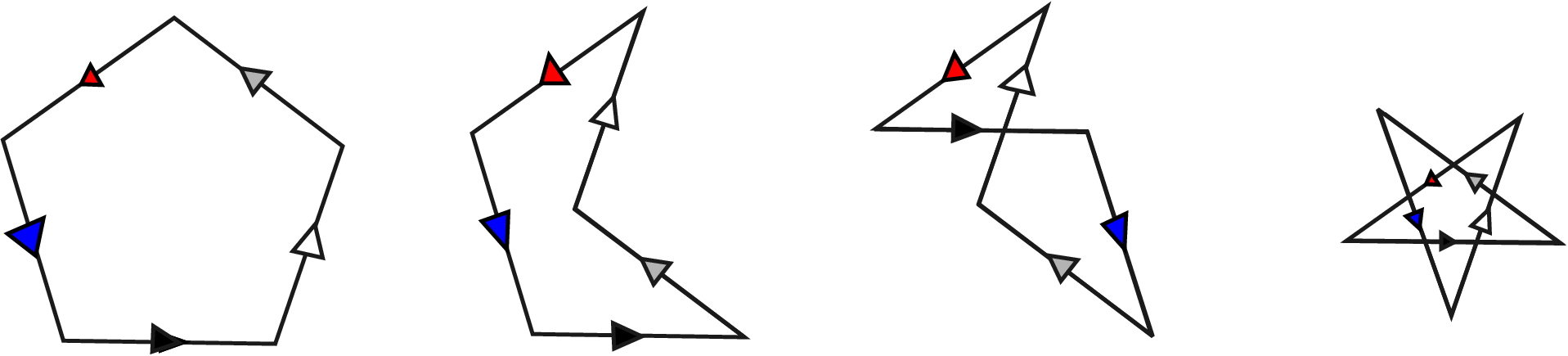}}
\newline
{\bf Figure 5.5:\/}  Regular pentagon permutations
\end{center}

We call an $N$-gon  {\it generic\/} if some
permutation makes it strictly convex. The set of
generic $N$-gons is open and dense, and invariant
under the action of $S_N$.
Topologically this subset is homeomorphic to
\begin{equation}
  f(N)=N!/(2N)
\end{equation}
copies of ${\cal C\/}_N$.   The reason why
$f(N)$ has the form it does is that the dihedral group,
which has order $2N$, permutes the convex classes.
We give a hyperbolic
structure to the subset of generic equilateral $N$-gons,
declaring it to be a disconnected union of $f(N)$
copies of ${\cal C\/}_N$.   By construction, $S_N$ acts
isometrically on our big space.  We call the
component corresponding to the strictly convex $N$-gons
the {\it convex component\/}.

We now enlarge
our space by taking the closures of all the components.
At the moment we still have a disjoint union of
hyperbolic polyhedra whose union (redundantly)
parametrizes all the equilateral $N$-gons.
Finally, we form an identification space by identifying
points in our union which represent the same
(equivalence class of) $N$-gon.  This is our hyperbolic structure on
the moduli space of equilateral $N$-gons.
We denote it by ${\cal A\/}_N$.
\newline
\newline
{\bf Remark:\/}
Technically, when $N$ is even, we have to
add to ${\cal A\/}_N$ the ideal points.  These
correspond to $N$-gons which lie in a single line.
For $N=6$ there are $10$ such.

\subsection{Pentagons}

In this section we explore
${\cal A\/}_5$ and prove Theorem \ref{penta}.
Before taking the quotient, we have
$12$ disjoint copies of ${\cal C\/}_5$, which is
a regular right angled hyperbolic pentagon.  These
pentagons are then glued edge-to-edge.
It turns out that $4$ are glued around each
vertex.  Figure 5.6 illustrates this for one
of the vertices of $C$, the copy of
${\cal C\/}_5$ that corresponds to
the convex pentagonal linkages. 
(I use the word {\it linkage\/} in this section to
avoid confusion;
the moduli space is also composed of pentagons.)
The vertices of
$C$ are certain isosceles triangles, in which
two pairs of consecutive edges point in the
same direction.
       
\begin{center}
\resizebox{!}{2.2in}{\includegraphics{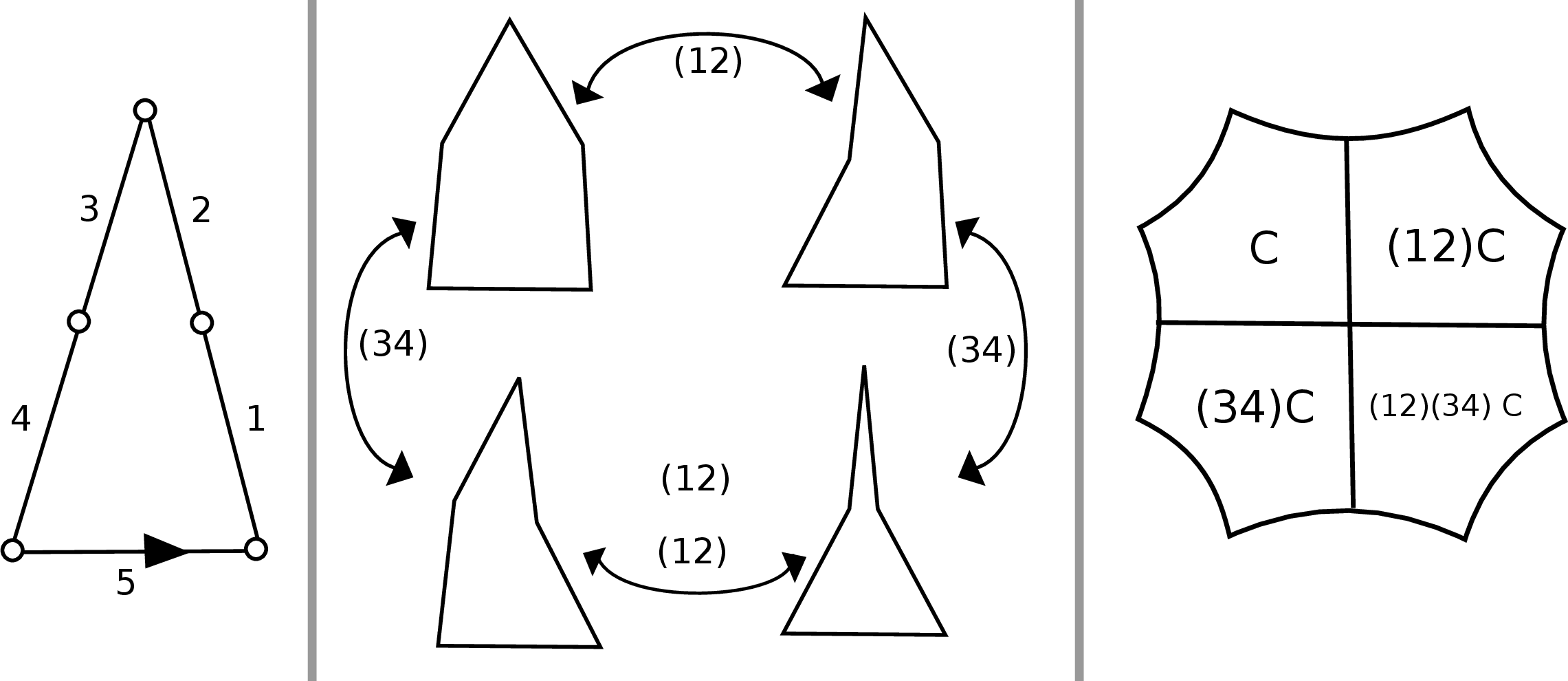}}
\newline
{\bf Figure 5.6:\/}  The local picture around a vertex.
\end{center}

The left panel of Figure 5.6 shows the vertex linkage.
The middle  panel shows the $4$ kinds of linkages which
contribute to the pentagons which glue together
around the vertex.  The right panel shows a hand-drawn
approximation of how the corresponding $4$ pentagons
would fit together if developed into the hyperbolic plane.
The subgroup generated by
the transpositions $(12)$ and $(34)$ fixes the
vertex and permutes the $4$ pentagons around it.
What makes this work is that $(12)$ and $(34)$
generate the dihedral subgroup of order $4$.

By symmetry, the picture is the same at every
vertex of our space.  Thus, the space we
get has a global hypebolic structure: It is isometric
to a very symmetric hyperbolic surface $\Sigma$.
The surface $\sigma$ has $12$ faces, and
$12 \times (5/2)=30$ edges and
$12 \times (5/4)=15$ vertices.  Hence
the surface has Euler characteristic
$\chi(\Sigma)=-3$.  
Given the classification of surfaces,
we can identify $\Sigma$ topologically as
the connected sum of a genus $2$ surface and
a projective plane.

\subsection{Hexagons}
\label{hexa}

In this section we explore ${\cal A\/}_6$ and
prove Theorem \ref{hex}.

Figure 5.7 shows the $5$ hexagons which
correspond to the $5$ vertices of ${\cal C\/}_6$.
The triangles correspond to vertices n
$\H^3$ and the segments correspond to
ideal vertices.

\begin{center}
\resizebox{!}{1.2in}{\includegraphics{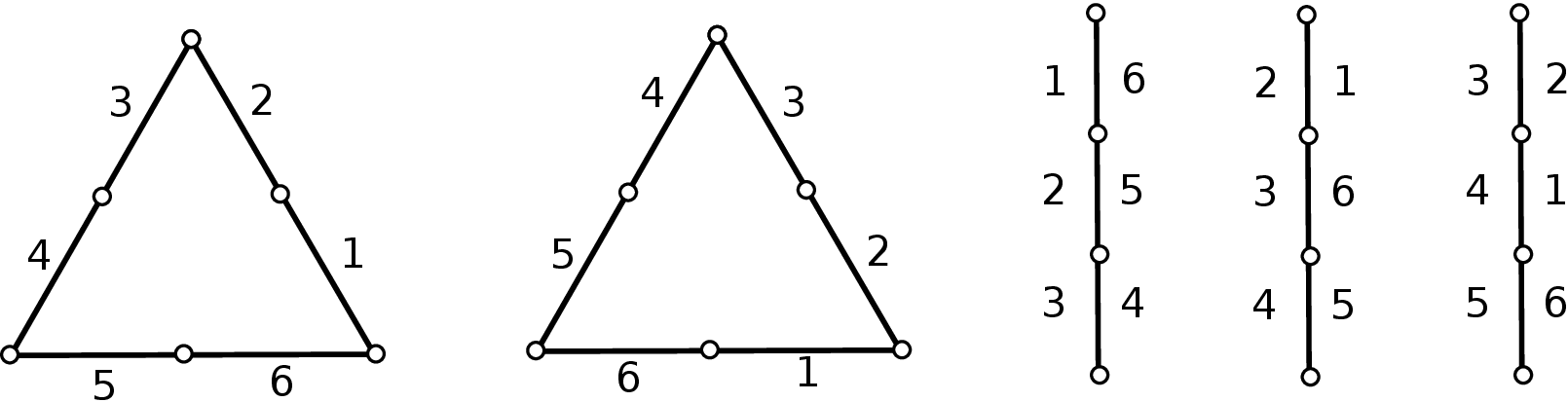}}
\newline
{\bf Figure 5.7:\/}  The vertics of ${\cal C\/}_6$.
\end{center}

The leftmost vertex is stabilized by the
order $8$ group $(\Z/2)^3$ generated by
the permutations $(12)$ and $(34)$ and $(56)$.
Thus, $8$ copies of ${\cal C\/}_6$ fit around
this vertex. Given the right-angles involved, the
identification space ${\cal A\/}_6$
is locally isometric to
$\H^3$ even along the vertices and edges.
The second triangular vertex has the same
kind of story, except that now the
permutations involves are $(23)$ and $(45)$ and $(61)$.

Figure 5.8 shows the hexagons corresponding to points
along the edge of
${\cal C\/}_6$ which connects the first two of the
ideal vertices shown above.

\begin{center}
\resizebox{!}{1.7in}{\includegraphics{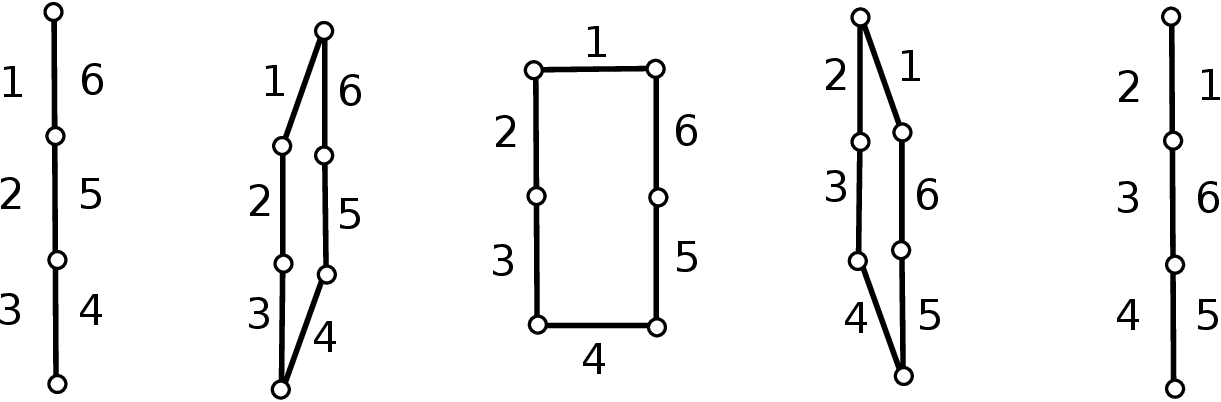}}
\newline
{\bf Figure 5.8:\/}  An edge of
${\cal C\/}_6$ connecting two ideal vertices.
\end{center}

This edge is stabilized by the order $4$ group
generated by $(23)$ and $(56)$.   From this, we
see that $4$ copies of ${\cal C\/}_6$ fit around
this edge.  Once again, given the right-angled
property of the faces of ${\cal C\/}_6$, this
means that our space ${\cal A\/}_6$ is locally
isometric to $\H^3$ around this edge.
The same story goes for the other two
edges of ${\cal C\/}_6$ that connect
ideal vertices.

Our analysis shows that ${\cal A\/}_6$ is locally isometric
to $\H^3$ in a neighborhood of every point of
${\cal C\/}_6$.   By symmetry, the same statement
holds for all points of ${\cal A\/}_6$.  Hence
${\cal A\/}_6$ is a hyperbolic $3$-manifold.
We can cut each copy of ${\cal C\/}_6$ in half
along the ideal triangle that is the convex hull of
the ideal vertices.   Each half is a pyramid obtained
by coning a regular ideal octahedron to the center
of mass.  In ${\cal A\/}_6$ we have $8$ of these
pyramids fitting together around each finite vertex to make
an ideal octahedron.
Thus ${\cal A\/}_6$ is tiled by regular ideal octahedra.
How many?

Well, ${\cal A\/}_6$ is obtained by gluing together
$f(6)=60$ copies of ${\cal C\/}_6$.
Each copy supplies $2$ pyramids, and we need
$8$ pyramids to make an ideal octahedron.
Thus, each copy of ${\cal C\/}_6$ supplies
$1/4$ of an octahedron.   We conclude
that ${\cal A\/}_6$ is tiled by $15$ regular
ideal hyperbolic octahedra.

\begin{center}
\resizebox{!}{1.5in}{\includegraphics{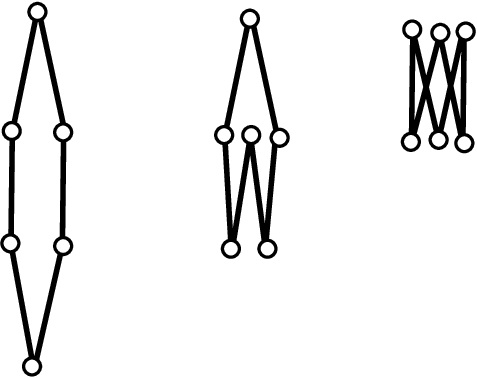}}
\newline
{\bf Figure 5.9:\/}  Shapes of nearly degenerate hexagons
\end{center}

The cusps of ${\cal A\/}_6$ are the degenerate hexagons
corresponding to the permutations of the
last 3 shown in Figure 5.7.
Figure 5.9 shows the representative shapes
of all $10$ degenerate hexagons.  There are
$3$ of the first kind, $6$ of the second kind,
and one of the third kind.  To get a
comprehensible picture we have
taken nearly degenerate hexagons rather
than actually degenerate ones.

Finally, Figure 5.10 shows the shapes of
the $15$ hexagons corresponding to the
centers of the ideal octahedra.  There are
respectively $2,6,3,3,1$ of these.

\begin{center}
\resizebox{!}{1in}{\includegraphics{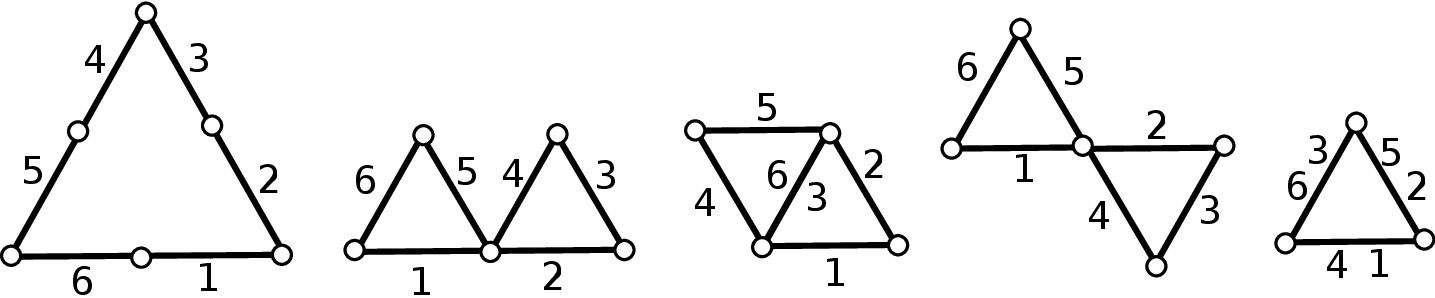}}
\newline
{\bf Figure 5.10:\/}  Shapes of the central hexagons
\end{center}

We have picked representative labelings.

\newpage

\section{References}

\noindent
[{\bf ABSST\/}]  P. Albers, G. Banhatti, P. Sadlo, R. Schwartz,
S. Tabachnikov, {\it Polygonal Symplectic Billiards\/},  J.
Experimental Math. (2025) to appear
\vskip 9 pt
\noindent
[{\bf ALW\/}] P. Albers, F. Lander, J. Westermann,
{\it Symplectic billiards for pair of polygons\/},
arXiv 2402.12244 (2024)
\vskip 9 pt
\noindent
[{\bf BDFI\/}] P Baird-Smith, D. Davis, E. Fromm, S Iyer,
{\it Tiling Billiards on Triangle Tiliings, and Interval Exchange
  Transformations\/},
Bulletin of the London Math Societ, Vol. 109, Issue 1 (2024)
\vskip 9 pt
\noindent
[{\bf BG\/}] C. Bavard and E. Ghys, {\it Polygones du plan et
  polyedres hyperboliques\/}, Geom. Dedicata 43(2),  (1992) pp 207--224
\vskip 9 pt
\noindent
[{\bf Cal\/}] D. Calegari, {\it Pentagonum Pentagonorum\/}, A.M.S.
Short Stories, Notices of the American Math Society, Vol. 68, No. 8
(2021)
\vskip 9 pt
\noindent
[{\bf KM\/}] M. Kapovich and J. Millson, {\it On the Moduli Space of
  Polygons in the Euclidean Plane\/},  J Diff. Geom, Vol. 42, No. 1,
(1995) 
\vskip 9 pt
\noindent
[{\bf R-G\/}] J. Richter-Gebert, {\it Mercator Workshop Lecture\/},
Heidelberg, June 2023
\vskip 9 pt
\noindent
[{\bf Sch1\/}] R. Schwartz, {\it Notes on Thurston's Shapes of
  Polyhedra\/},  \newline
arXiv 1506.07252 (2015)
\vskip 9 pt
\noindent
[{\bf Sch2\/}] R. Schwartz, {\it The Pentagram Map\/}
Exp. Math., Vol 1 (1992) pp 85--90
\vskip 9 pt
\noindent
[{\bf Si\/}] J. Silverman, private communication, May 1, 2025.
\vskip 9 pt
\noindent
[{\bf Sp\/}], B. Springborn,
{\it A unique representation of polyhedral typs.
  Centering via M\"obius transformations\/},
Math. Zeit. {\bf 249\/} (2004) pp 513-517
\vskip 9 pt
\noindent
[{\bf St\/}]K. Stephenson, {\it Circle Packing: A Mathematical
  Tale\/},
Bulletin of the A.M.S., Vol. 50, No. 11 (2003) pp 1376 -- 1387
\vskip 9 pt
\noindent
[{\bf T\/}] W. Thurston, {\it Shapes of Polyhedra and triangulations
  of the sphere\/}, \newline
Geometry $\&$ Topology Monographs, Vol. 1:
The Epstein birthday shrift, pp 511-549  (1998)  See also
arXiv:math/9801088.
\vskip 9 pt
\noindent
[{\bf W\/}] S. Wolfram et. al., {\it Mathematica\/}, Vers. 11,
Wolfram Res. Inc. (2024)

\end{document}